\documentclass[12pt,twoside,a4paper]{article}
\usepackage[symbol]{footmisc}
\renewcommand{\thefootnote}{\arabic{footnote}}%
\usepackage{bm}
\usepackage{braket}
\usepackage{graphicx}
\usepackage{graphics}
\usepackage{theorem}
\usepackage{empheq}
\usepackage{amsmath}
\usepackage{amssymb,mathrsfs}
\usepackage{latexsym}
\usepackage[subrefformat=parens]{subcaption}
\usepackage{cases}
\usepackage{placeins} %To control the float position
\usepackage[round,authoryear]{natbib}	
\usepackage{bbm}   	%\mathbbm{12}

\newcommand{\one}{\mathbbm{1}}
\usepackage{color}	
\usepackage{algorithm}
\usepackage{algpseudocode}
\usepackage{caption}
\newfloat{procedure}{H}{loa}
\floatname{procedure}{Procedure} %It may be done better
\theorembodyfont{\itshape}
\newtheorem{thm}{Theorem}[section]
\newtheorem{lem}[thm]{Lemma}%[section]
\newtheorem{prop}[thm]{Proposition}%[section]
\newtheorem{coro}[thm]{Corollary}%[section]
\newtheorem{Power-law}[thm]{Power-law Formulas}%[section]

\theorembodyfont{\rmfamily}
\newtheorem{defn}[thm]{Definition}%[section]
\newtheorem{assumpt}[thm]{Assumption}%[section]

\newtheorem{rem}[thm]{Remark}%[section]
\newcommand{\qed}{\hspace*{\fill}$\Box$}
\numberwithin{equation}{section}  % If you number theorems, etc. within sections,
\makeatletter
\def\eqnarray{\stepcounter{equation}\let\@currentlabel=\theequation
\global\@eqnswtrue
\global\@eqcnt\z@\tabskip\@centering\let\\=\@eqncr
$$\halign to \displaywidth\bgroup\@eqnsel\hskip\@centering
  $\displaystyle\tabskip\z@{##}$&\global\@eqcnt\@ne 
  \hfil$\;{##}\;$\hfil
  &\global\@eqcnt\tw@ $\displaystyle\tabskip\z@{##}$\hfil 
   \tabskip\@centering&\llap{##}\tabskip\z@\cr}
\makeatother
\makeatletter
\newif\if@borderstar
\def\bordermatrix{\@ifnextchar*{%
 \@borderstartrue\@bordermatrix@i}{\@borderstarfalse\@bordermatrix@i*}%
}
\def\@bordermatrix@i*{\@ifnextchar[{\@bordermatrix@ii}{\@bordermatrix@ii[()]}}
\def\@bordermatrix@ii[#1]#2{%
\begingroup
 \m@th\@tempdima8.75\p@\setbox\z@\vbox{%
 \def\cr{\crcr\noalign{\kern 2\p@\global\let\cr\endline }}%
 \ialign {$##$\hfil\kern 2\p@\kern\@tempdima & \thinspace %
  \hfil $##$\hfil && \quad\hfil $##$\hfil\crcr\omit\strut %
  \hfil\crcr\noalign{\kern -\baselineskip}#2\crcr\omit %
  \strut\cr}}%
 \setbox\tw@\vbox{\unvcopy\z@\global\setbox\@ne\lastbox}%
 \setbox\tw@\hbox{\unhbox\@ne\unskip\global\setbox\@ne\lastbox}%
 \setbox\tw@\hbox{%
  $\kern\wd\@ne\kern -\@tempdima\left\@firstoftwo#1%
  \if@borderstar\kern 2pt\else\kern -\wd\@ne\fi%
 \global\setbox\@ne\vbox{\box\@ne\if@borderstar\else\kern 2\p@\fi}%
 \vcenter{\if@borderstar\else\kern -\ht\@ne\fi%
  \unvbox\z@\kern -\if@borderstar2\fi\baselineskip}%
 \if@borderstar\kern-2\@tempdima\kern2\p@\else\,\fi\right\@secondoftwo#1 $%
 }\null \;\vbox{\kern\ht\@ne\box\tw@}%
\endgroup
}
\makeatother
\makeatletter
\DeclareRobustCommand\widecheck[1]{{\mathpalette\@widecheck{#1}}}
\def\@widecheck#1#2{%
    \setbox\z@\hbox{\m@th$#1#2$}%
    \setbox\tw@\hbox{\m@th$#1%
       \widehat{%
          \vrule\@width\z@\@height\ht\z@
          \vrule\@height\z@\@width\wd\z@}$}%
    \dp\tw@-\ht\z@
    \@tempdima\ht\z@ \advance\@tempdima2\ht\tw@ \divide\@tempdima\thr@@
    \setbox\tw@\hbox{%
       \raise\@tempdima\hbox{\scalebox{1}[-1]{\lower\@tempdima\box
\tw@}}}%
    {\ooalign{\box\tw@ \cr \box\z@}}}
\makeatother
\newcommand{\ol}{\overline}

\newcommand{\vc}{\bm}
\newcommand{\vmax}{\vee}

\newcommand{\dm}{\displaystyle}

\newcommand{\varep}{\varepsilon}

\newcommand{\EE}{\mathsf{E}}
\newcommand{\PP}{\mathsf{P}}

\newcommand{\calR}{\mathcal{R}}

\newcommand{\bb}{\mathbb}

\newcommand{\bbN}{\mathbb{N}}
\newcommand{\bbR}{\mathbb{R}}
\newcommand{\bbS}{\mathbb{S}}

\renewcommand{\thefootnote}{\alph{footnote}}
\makeatother

\renewcommand{\thefootnote}{\fnsymbol{footnote}}
\begin{document}\thispagestyle{empty} 

\hfill
%{\small Last update date: \today}
%Submitted to STOCHASTIC MODELS, \today.

\vspace{-10mm}

{\large{\bf
\begin{center}
Curse of scale-freeness: Intractability of large-scale optimization with multi-start methods%
\footnote{Published in \textit{Extremes} on 21 April 2025. DOI: https://doi.org/10.1007/s10687-025-00511-w}
\end{center}
}
}

\setcounter{footnote}{0}
\renewcommand{\thefootnote}{\arabic{footnote}}
%%%%%%%%%%%%%%%%%%%%%%%%%%%%%%%%%%%%%%%%%%%%%%%%%%%%%%%%%%%%%%%%%
\begin{center}
{
Hiroyuki Masuyama\footnote[1]{E-mail: masuyama@tmu.ac.jp}
}

{\small \it
Graduate School of Management, Tokyo Metropolitan University\\
Tokyo 192--0397, Japan
}

\bigskip
%%%%%%%%%%%%%%%%%%%%%%%%%%%%%%%%%%%%%%%%%%%%%%%%%%%%%%%%%%%%%%%%%
{
Hiroshige Dan\footnote[2]{E-mail: dan@waseda.jp}
}

{\small \it
School of Creative Science and Engineering, Shinjuku-ku, Tokyo, 169--8555, Japan
}% \small ends

%%%%%%%%%%%%%%%%%%%%%%%%%%%%%%%%%%%%%%%%%%%%%%%%%%%%%%%%%%%%%%%%%
\bigskip
{
Shunji Umetani%
\footnote[3]{E-mail: umetani@ist.osaka-u.ac.jp}
}

{\small\it
Department of Information and Physical Sciences, Graduate School of Information Science and Technology, The University of Osaka, Suita, Osaka, 565--0871, Japan
}%
%%%%%%%%%%%%%%%%%%%%%%%%%%%%%%%%%%%%%%%%%%%%%%%%%%%%%%%%%%%%%%%%%
\setcounter{footnote}{3}

\bigskip
\medskip

{\small
\textbf{Abstract}

\medskip

\begin{tabular}{p{0.85\textwidth}}
This paper investigates the intractability of large-scale optimization with multi-start methods. For the theoretical performance analysis, we focus on random multi-start (RMS), which is one of the representative multi-start methods, including RMS local search and greedy randomized adaptive search procedure (GRASP). Our primary theoretical contribution is to derive, by using extreme value theory, power-law formulas for the two quantities: (i) the expected improvement rate of the best empirical objective value (EOV); (ii) the expected relative gap between the best EOV and the supremum of the EOVs. Notably, the expected relative gap exhibits scale-freeness as a function of the number of iterations. Consequently, the half-life of the expected relative gap is asymptotically proportional to the number of iterations executed by the RMS method. This result can be interpreted as the {\it curse of scale-freeness}---a Zeno's paradox-like phenomenon---expressed by the metaphor ``{\it Reaching for the goal makes it slip away.}" Through numerical experiments, we observe that several RMS algorithms applied to traveling salesman problems suffer from the curse of scale-freeness. Furthermore, we show that overcoming this curse requires a powerful local search algorithm with effective restart and diversification strategies that exponentially accelerate solution improvement relative to the RMS method.
\end{tabular}
}
% \samllsize ends
\end{center}

\begin{center}
\begin{tabular}{p{0.90\textwidth}}
{\small
{\bf Keywords:} %
multi-start method; extreme value theory (EVT); scale-free; power law; traveling salesman problem (TSP); Zeno's paradox
% 
% End of Keywords
%

\medskip

{\bf Mathematics Subject Classification:} %
90C27; 60G70; 62G32
}%\samllsize ends
\end{tabular}

\end{center}

%%%%%%%%%%%%%%%%%%%%%%%%%%%%%%%%%%%%%%%%%%%%%%%%%%%%%%%%%%%%%%%%%%%%%
\section{Introduction}\label{sec_introduction}

The purpose of this paper is to statistically reveal the intractability of large-scale optimization with the multi-start method \citep{Mart13,Mart18}. The multi-start method is a standard metaheuristic strategy that iteratively restarts a local search (LS) algorithm with a diversification strategy. In terms of diversification strategies, multi-start methods can be classified into two types: (i) random multi-start (RMS) methods (see Algorithm~\ref{algo_RMS}), such as random multi-start local search and greedy randomized adaptive search procedure (GRASP) \citep{Res16}; and (ii) random perturbation methods (see Algorithm~\ref{algo_RP}), such as iterated local search (ILS) \citep{Lou19}, variable neighborhood search \citep{Han01}, and memetic algorithms \citep{Neri12}. The RMS method generates an initial solution from scratch by using a randomized construction algorithm, whereas random perturbation methods do so by perturbing an existing {\it good} (often locally optimal) solution. In general, the multi-start method achieves relatively good performance despite its simple implementation, making it a standard benchmark for designing heuristic algorithms in many real-world applications.

\begin{algorithm}[h]
 \caption{Random multi-start (RMS) for maximization problems}\label{algo_RMS}
\begin{algorithmic}[1]
\State $t \gets 1$ \Comment{$t$ is the number of trials.}
\State $x^{\dag} \gets -\infty$
\Comment{$x^{\dag}$ is the best empirical objective value (EOV).}
\While{$t \le T$} \Comment{$T$ is the maximum number of trials.}
\State $\vc{s}^0 \gets \mathrm{Random}$ \Comment{Randomly generate an initial solution $\vc{s}^0$.}
\State $\vc{s} \gets \mathrm{LS}(\vc{s}^0)$ \Comment{Generate an empirical solution $\vc{s}$ from $\vc{s}^0$ by LS.}
\If {$x(\vc{s}) > x^{\dag}$} \Comment{$x(\vc{s})$ is the objective value of the empirical solution $\vc{s}$.}
\State $\vc{s}^{\dag} \gets \vc{s}$ 
\Comment{Update the best empirical solution $\vc{s}^{\dag}$ with $\vc{s}$.}
\State $x^{\dag} \gets x(\vc{s}^{\dag})$  
\Comment{Update the best EOV $x^{\dag}$ with $x(\vc{s}^{\dag})$.}
\EndIf
\State $t \gets t+1$
\EndWhile
\State \Return $\vc{s}^{\dag}$
\end{algorithmic}
\end{algorithm}

\begin{algorithm}[h]
 \caption{Random perturbation for maximization problems}\label{algo_RP}
\begin{algorithmic}[1]
\State $t \gets 1$ 
\State $\vc{s}^0 \gets \mathrm{Initialize}$ \Comment{Randomly generate an initial solution $\vc{s}^0$ for the first trial.}
\State $\vc{s}^{\dag} \gets \mathrm{LS}(\vc{s}^0)$
\Comment{Generate an empirical solution $\vc{s}^{\dag}$ from $\vc{s}^0$ by LS.}
\State $x^{\dag} \gets x(\vc{s}^{\dag})$
\Comment{Set $x^{\dag}= x(\vc{s}^{\dag})$ as the first (and therefore best) EOV.}
\While{$t \le T$} 
\State $\vc{s}^0 \gets \mathrm{Kick}(\vc{s}^{\dag})$ 
\Comment{Generate an initial solution $\vc{s}^0$ by randomly perturbing $\vc{s}^{\dag}$.}
\State $\vc{s} \gets \mathrm{LS}(\vc{s}^0)$ \Comment{Generate an empirical solution $\vc{s}$ from $\vc{s}^0$ by LS.}
\If {$x(\vc{s}) > x^{\dag}$}
\State $\vc{s}^{\dag} \gets \vc{s}$
\Comment{Update the best empirical solution $\vc{s}^{\dag}$ with $\vc{s}$.}
\State $x^{\dag} \gets x(\vc{s}^{\dag})$
\Comment{Update the best EOV $x^{\dag}$ with $x(\vc{s}^{\dag})$.}
\EndIf
\State $t \gets t+1$
\EndWhile
\State \Return $\vc{s}^{\dag}$
\end{algorithmic}
\end{algorithm}

This paper clarifies the intractability of large-scale optimization with the multi-start method in the following way. First, we apply extreme value theory (EVT) to the process of solving a large-scale optimization problem by the RMS method (Algorithm~\ref{algo_RMS}), which is straightforward to analyze mathematically. We then derive power-law formulas for the improvement process of the empirical objective values (EOVs)---where ``objective value" is an abbreviation for ``objective function value"---generated by iteratively running the RMS method. These power-law formulas provide a theoretical basis for understanding the intractability of large-scale optimization, which we figuratively express as the {\it curse of scale-freeness}. Furthermore, through numerical experiments---rather than theoretical analysis using EVT---we confirm that even the ILS method (Algorithm~\ref{algo_RP}), commonly regarded as a high-performance random perturbation method, cannot readily break the curse of scale-freeness.

\begin{figure}[h]%%%Figure 1%%%
\centering
\includegraphics [scale=0.45, bb=0 0 878 200]{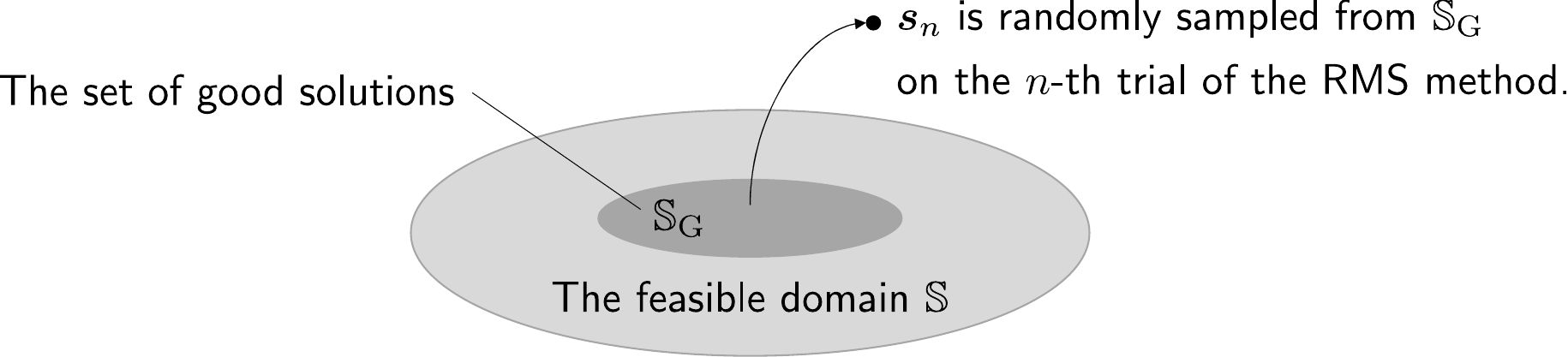}
\caption{The RMS method as a random sampler from the set $\bbS_{\rm G}$ of good solutions.}\label{fig_image_RMS}
\end{figure}

As in \cite{Gidd14}, this paper considers a class of large-scale mixed-integer programming problems described by a feasible domain $\bbS$, a set of decision variables $\vc{s} \in \bbS$, and a nonnegative objective function $x(\vc{s})$ to be maximized:
\begin{equation}\label{our_prob}
\left\{ \,
\begin{alignedat}{2}
\mbox{maximize}   &  & \quad	    x(\vc{s}) &\in \bbR_+:=[0,\infty)
\\
\mbox{subject to} &  & \quad  	\vc{s}    & \in \bbS.
\end{alignedat}
\right.
\end{equation}
The following assumptions are made throughout this paper.
\begin{assumpt}\label{assumpt_basic}
The RMS method generates a single empirical solution $\vc{s}_n \in \bbS$ of Problem (\ref{our_prob})\footnote{In the published version, ``Problem Eq.\ 1.1'' appears instead of the intended ``Problem 1.1'' (although ``Problem (1.1)'' would be more natural, as in this preprint). The correction, requested during proofreading, was not implemented.} for each trial $n \in \bbN:=\{1,2,\dots\}$, forming the set of {\it good} (typically locally optimal) solutions $\bbS_{\rm G} = \{\vc{s}_n; n\in \bbN\}$ (see Figure~\ref{fig_image_RMS}). We assume that $\bbS$ and its subset $\bbS_{\rm G}$ are infinite or sufficiently large to be considered infinite. Thus, the EOVs $X_n:=x(\vc{s}_n)$, $n \in \bbN$, are independent and identically distributed (i.i.d.), implying that the RMS method is a random sampler from $\bbS_{\rm G}$ (\citealt{Gold77,Gold78,Gold79}). Let $X$ denote a generic random variable from the i.i.d.\ sequence $\{X_n\}$ with distribution function $F$:
\begin{align*}
F(x) = \PP(X \le x) = \PP(X_n \le x),\qquad x \in \bbR_+.
%\label{ean_F(x)}
\end{align*}
Let $x^*$ denote the right endpoint of $F$, i.e.,
\begin{align}
x^* = \sup\{x \in \bbR_+: F(x) < 1\},
\label{defn_x^*}
\end{align}
which is the supremum of the EOVs and is assumed to be positive.
\end{assumpt}

\begin{rem}\label{rem_SEV}
If the feasible domain $\bbS$ is finite and the employed RMS algorithm is powerful enough to reach the (global) optimal solution given a suitable initial solution, then the supremum $x^*$ of the EOVs coincides with the optimal value. Otherwise, $x^*$ may be infinite in some cases, or finite yet not equal to the optimal value (the latter can occur when $\bbS$ is open).
\end{rem}

\medskip

\begin{rem}
Usually, each time the RMS method reaches a locally optimal solution, it restarts the local search with an initial solution generated in the same stochastic manner, independent of the computational history. Consequently, the i.i.d.\ property of the EOVs generated by the RMS method is a reasonable assumption, given a sufficiently large number of possible initial solutions.
\end{rem}

\medskip

A topic related to our research is the estimation of the optimal value, which may or may not coincide with the supremum $x^*$ of the EOVs (see Remark~\ref{rem_SEV}). Estimation methods for the optimal value can typically be classified into two types. The first is the relaxation method (see, e.g., \citealt{Held70,Fish75}), which solves relaxation problems to establish bounds on the optimal value. The second is the statistical method, which employs statistical tools, such as EVT, to provide point and interval estimates of the optimal value. Pioneering work on the EVT-based estimation of the optimal value was conducted by \cite{Gold77} (see also \citealt{Gold78,Gold79}). \cite{Gidd14} provides a comprehensive survey not only of EVT-based estimation methods but also of other statistical ones. Following \cite{Gidd14}, \cite{Carl15,Carl16} evaluated the quality of heuristic solutions by estimating the optimal value with the Weibull distribution. Several other papers have discussed the EVT-based optimal value estimation for specific optimization problems; see Section 1 of \cite{Vela22}. However, \cite{Vela22} raises concerns about the reliability of such EVT-based estimation methods.

\begin{figure*}[htb]%%%Figure 2%%%
\begin{tabular}{cc}
\begin{minipage}[t]{0.45\hsize}
\centering
\includegraphics [scale=0.45,bb=0 0 461 346]{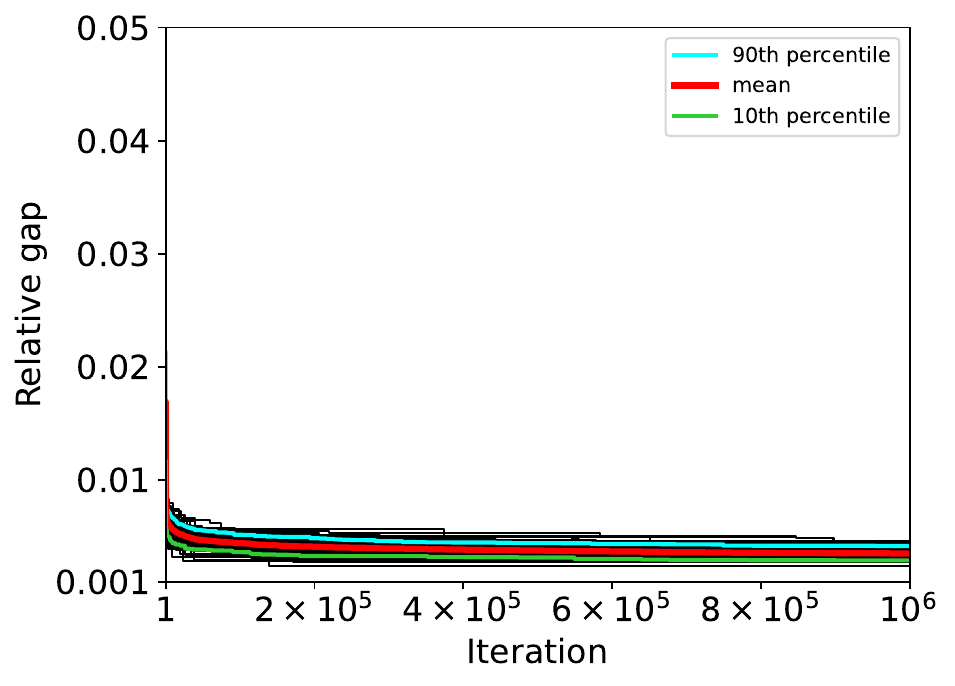}
\subcaption{Linear axes}
\label{fig_RE-NN-LKH-normal}
\end{minipage} &
\quad
\begin{minipage}[t]{0.45\hsize}
\centering
\includegraphics [scale=0.45,bb=0 0 461 346]{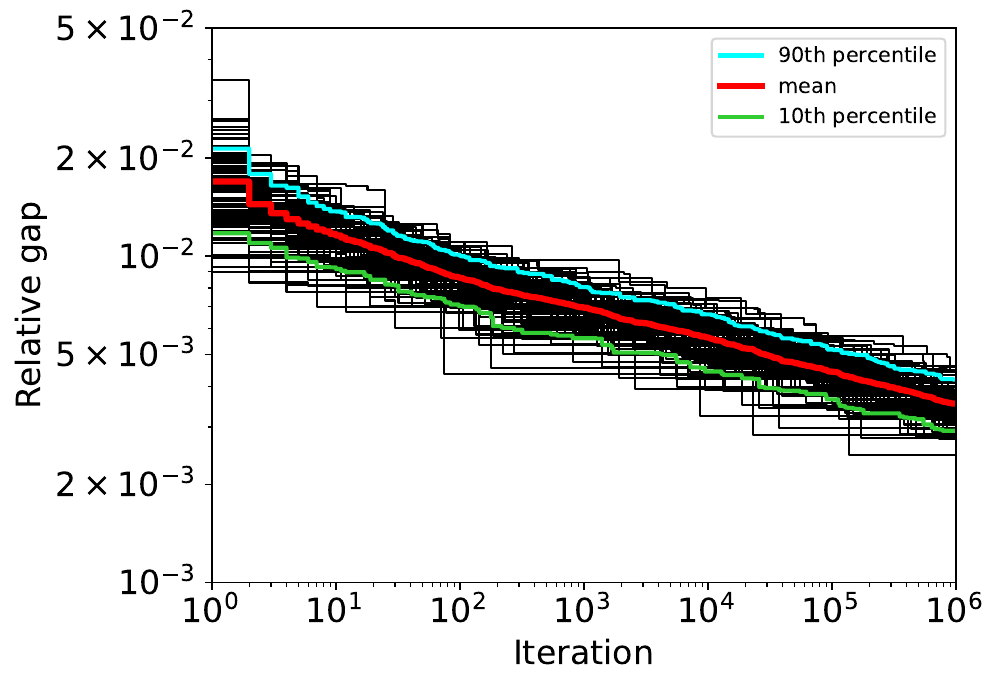}
\subcaption{Logarithmic axes}
\label{fig_RE-NN-LKH-log}
\end{minipage}
\end{tabular}
\caption{
Evolution of the relative gap between the best EOV and the optimal value for an RMS algorithm applied to 100 random TSP (Traveling Salesman Problem) instances. Figures~\ref{fig_RE-NN-LKH-normal} and \ref{fig_RE-NN-LKH-log} show the same data, but Figure~\ref{fig_RE-NN-LKH-normal} uses linear axes, whereas Figure~\ref{fig_RE-NN-LKH-log} uses logarithmic axes.}\label{fig_power-law}
\end{figure*}

Accurate estimation of the gap between the best EOV and the optimal value is insufficient for predicting the computational effort needed to close the gap. In many cases, this gap narrows rapidly during the early stages of computation but soon reaches a {\it steady state} where further reduction rarely occurs. Figure~\ref{fig_RE-NN-LKH-normal} illustrates a typical example of this behavior.

Figure~\ref{fig_RE-NN-LKH-log} (the log-log plot of Figure~\ref{fig_RE-NN-LKH-normal}) shows the approximately linear evolution of the relative gap (error) of the best EOV, measured from the optimal value. Thus, the power-law behavior---or equivalently, scale-freeness---is evident in this evolution. In other words, the number of additional trials required to halve the relative gap increases asymptotically in proportion to the number of trials completed. This Zeno's paradox-like phenomenon can be figuratively described as the {\it curse of scale-freeness: ``Reaching for the goal makes it slip away."} As discussed in detail later, the curse of scale-freeness represents the intractability of large-scale optimization, at least when using the RMS method.

The above considerations are supported by the power-law formulas for the evolution of the relative gap of the best EOV, which are summarized below.

\smallskip

\noindent
{\bf Summary of the power-law formulas:~}
Let $Z_n = \max(X_1,X_2,\dots,X_n)$ for $n \in \bbN$, and let $\vmax_{i=k}^{\ell}X_i=\dm\max(X_k,X_{k+1},\dots,X_{\ell})$ for $k,\ell \in \bbN$ with $k < \ell$. Since $\{X_n\}$ is i.i.d., the random variable $\vmax_{i=N+1}^{N+n}X_i$ is independent of $Z_N$ and has the same distribution as $Z_n$. Furthermore, let
\begin{alignat}{2}
R_n(x) &= {(Z_n - x)_+ \over x}, & \quad  0 &< x < x^*,\ n \in \bbN,
\label{defn_R_n(x)}
\end{alignat}
where $(x)_+ = \max(x,0)$ for $x \in \bbR:=(-\infty,\infty)$. Let
\begin{alignat}{2}
\varDelta_n(x)
&=  {x^* - \max(Z_n,x) \over x^*}
= {x^* - x \over x^*} - {(Z_n-x)_+ \over x^*}, & \quad  0 &< x < x^*,\ n \in \bbN,
\label{defn_Delta_n(x)}
\end{alignat}
where $x^* < \infty$ is always assumed when considering $\varDelta_n(x)$. It follows from (\ref{defn_R_n(x)}) and (\ref{defn_Delta_n(x)}) that $R_n(x)$ and $\varDelta_n(x)$ are related by the affine equation:
\begin{align}
\varDelta_n(x) 
= {x^* - x \over x^*} - { x \over x^* } R_n(x)
= {x \over x^*} \left( {x^* - x \over x} - R_n(x) \right).
\label{eqn_affine_relation}
\end{align}
Note that, given that the best EOV equals $x$, (i) $\EE[R_n(x)]$ represents the (conditional) expected improvement rate of the best EOV after $n$ additional trials of the RMS method, and (ii) $\EE[\varDelta_n(x)]$ represents the (conditional) expected relative gap of the best EOV after $n$ additional trials, measured from the supremum $x^* < \infty$ of the EOVs. According to Corollary~\ref{coro_regularly-varying}, the two expectations, $\EE[R_n(x)]$ and $\EE[\varDelta_n(x)]$, satisfy the power-law formulas given below:
\begin{subequations}\label{power_law}
\begin{alignat}{2}
\EE[R_n(x)] &= {x^* - x \over x} - {x^*  \over x} \EE[\varDelta_n(x)],
& \qquad 0 &< x < x^*, 
\label{power_law_01}
\\
\EE[\varDelta_n(x)] &= L(n) n^{\xi},
& \qquad 0 &< x < x^*, 
\label{power_law_02}
\end{alignat}
\end{subequations}
where $\xi < 0$ is a constant (in fact, (\ref{power_law}) also holds when $\xi = 0$; see Corollary~\ref{coro_regularly-varying}) and $L$ is a slowly varying function (see Definition~\ref{defn_class_R}). The above results are presented in detail in Section~\ref{sec_main}, along with other related results.

\smallskip

The power-law formulas in (\ref{power_law}) with a negative exponent $\xi < 0$ provide the theoretical foundation for the curse of scale-freeness. Equation~(\ref{power_law_01}) shows that additional trials increase $\EE[R_n(x)]$ toward its supremum $(x^* - x)/x$ at the same rate as the decrease in $(x^*/x)\EE[\varDelta_n(x)]$, which decays toward zero at a power-law rate of order $n^{\xi}$. Notably, (\ref{power_law_02}) shows that $\EE[\varDelta_n(x)]$ is asymptotically independent of $x$ and exhibits a scale-free property: For $x \in (0,x^*)$,
\begin{align*}
\lim_{n\to\infty} {\EE[\varDelta_{cn}(x)] \over \EE[\varDelta_n(x)]}
&=  c^{\xi}, \qquad 0 < x < x^*,~c \in \bbN,
\end{align*}
which yields
\begin{align}
 \lim_{n\to\infty} 
{\EE[\varDelta_{\lceil 2^{-1/\xi} \rceil n}(x)] \over \EE[\varDelta_n(x)]}
=\lceil 2^{-1/\xi} \rceil^{\xi}
\le
{1 \over 2}
\le
\lfloor 2^{-1/\xi} \rfloor^{\xi}
=
\lim_{n\to\infty} 
{\EE[\varDelta_{\lfloor 2^{-1/\xi} \rfloor n}(x)] \over \EE[\varDelta_n(x)]},
\label{eqn_curse_of_scale_freenss}
\end{align}
where $\lfloor \, \cdot \, \rfloor$ and $\lceil \, \cdot \, \rceil$ denote the floor and ceiling functions, respectively. As illustrated by Figure~\ref{fig_RE-NN-LKH-log}, Equation~(\ref{eqn_curse_of_scale_freenss}) implies that the {\it half-life} of the expected relative gap is asymptotically proportional to the number of trials completed. This phenomenon is aptly described as ``{\it Reaching for the goal makes it slip away.}"

Our power-law formulas (\ref{power_law_01}) and (\ref{power_law_02}) imply that overcoming the curse of scale-freeness depends on devising an LS algorithm that incorporates advanced restart and diversification strategies to achieve exponential acceleration in improving the best empirical solution relative to the RMS method (see Section~\ref{subsec_intractability}). These results indicate that the curse of scale-freeness is very difficult to overcome. To confirm this conclusion, we apply three RMS and three ILS algorithms to five instances of the Traveling Salesman Problem (TSP) from TSPLIB (see Table~\ref{TSPLIB_instances} and Appendix~\ref{sec_num_setup} for details) and summarize the numerical results in Figures~\ref{RMS-ILS_RA-LKH}--\ref{RMS-ILS_GR-LKH}. The figures reveal that the relative gap between the best EOV and the optimal value decays at most at a power-law rate, demonstrating that the ILS algorithms cannot overcome the curse of scale-freeness. Since the ILS method is known as one of the most powerful metaheuristics among multi-start methods, it is likely that other multi-start algorithms are also trapped by the curse of scale-freeness.

The rest of this paper is organized as follows. Section~\ref{sec_preliminaries} presents the preliminaries for our EVT analysis of empirical solutions generated by the RMS method. Section~\ref{sec_main} provides the main theoretical results of the paper, including the power-law formulas described above. Based on these theoretical results, Section~\ref{sec_discussion} discusses the intractability of large-scale optimization with multi-start methods. Finally, Section~\ref{sec_remarks} offers concluding remarks.

\section{Preliminaries}\label{sec_preliminaries}

This section is divided into two subsections. Section~\ref{subsec_basic} introduces the fundamental assumption of EVT, which forms the basis for the theoretical arguments presented later. Section~\ref{subsec_class_R} presents the preliminary results on regularly varying functions that are necessary for deriving the power-law formulas for the expected improvement rate and the expected relative gap discussed in the next section.

For later reference, we introduce some symbols for functions and sets of numbers. For any eventually positive function $f$ on $\bbR$, the notation $g(x)=o(f(x))$ means that
\begin{align*}
\lim_{x\to\infty}{|g(x)| \over f(x)} = 0.
\end{align*}
In addition, for any non-decreasing function $f$, let $f^{\gets}$ denote the {\it left-continuous inverse} of $f$, i.e.,
\begin{align*}
f^{\gets}(x) = \inf\{y \in \bbR: f(y) \ge x\}.
\end{align*}
Finally, let $\one(\,\cdot\,)$ denote the indicator function, which equals one if the condition in the parentheses is true, and equals zero otherwise.

\subsection{The fundamental assumption of EVT}\label{subsec_basic}

This subsection introduces the fundamental assumption of EVT \cite[Section~1.1]{Haan06} in addition to Assumption~\ref{assumpt_basic}. The fundamental assumption is reformulated as two equivalent conditions, which form the basis for deriving the theoretical results presented in the following sections.

The fundamental assumption of EVT is as follows:
\begin{assumpt}\label{assumpt_GEV}
There exist two sequences of constants $\{ a_n>0;n\in \bb{N}\}$ and $\{ b_n\in\bbR;n\in \bb{N}\}$ and a non-degenerate distribution $G$ such that 
\begin{align*}
\lim_{n\to\infty}
\PP\!\left( {Z_n - b_n \over a_n} \le z \right) 
&= \lim_{n\to\infty} [F(a_n z + b_n)]^n
= G(z),
\end{align*}
for every continuity point $z$ of $G$. In this case, $F$ is said to be in the maximum domain of attraction (MDA) of $G$, denoted by $F \in \mathsf{MDA}(G)$.
\end{assumpt}

\medskip

For later discussion, we present two conditions that are equivalent to Assumption~\ref{assumpt_GEV}. To describe these conditions, we introduce a specific function $U$ and the standard GEV distribution as follows. Let $U$ denote a function such that
\begin{align}
U(t) 
&= \left( {1 \over 1-F} \right)^{\gets}\!\!(t)
= \inf\left\{x \in \bbR: {1 \over 1 - F(x)} \ge t \right\},\qquad t \ge 1.
\label{defn_U}
\end{align}
Equations (\ref{defn_U}) and (\ref{defn_x^*}) imply that
\begin{align}
\lim_{t\to\infty} U(t) = x^*.
\label{eqn_lim_U(t)=x^*}
\end{align}
Furthermore, let $G_{\xi}$ denote the {\it standard generalized extreme value (GEV) distribution}, defined as follows:
\begin{equation}
G_{\xi}(z)
=
\exp\left\{ 
-( 1 + \xi z)_+^{-1/\xi} \right\},
\qquad z \in \bbR,
\label{defn_G_{xi}}
\end{equation}
where $(1 + \xi z)^{-1/\xi}$ is interpreted as $e^{-z}$ if $\xi=0$ (this convention comes from \citealt[Theorem~1.1.3]{Haan06}) and thus 
\begin{align*}
G_0(z) = \exp\{-e^{-z} \}, \qquad z \in \bbR.
\end{align*}
Note that the mean $\ol{m}_{\xi} := \int_{-\infty}^{\infty} z\,  G_{\xi}'(z) \,d z$ is given by
\begin{subequations}\label{eqn_ol{m}_{xi}}
\begin{empheq}[left = {\ol{m}_{\xi} = \empheqlbrace \,}]{alignat = 2}
& -\varGamma(-\xi) - {1 \over \xi}, &\qquad &   \mbox{$\xi \neq 0,~\xi < 1$},  
\label{eqn_ol{m}_{xi}-a} 
\\
& \gamma,                           &\qquad &   \mbox{$\xi = 0$},  
\label{eqn_ol{m}_{xi}-b}
\end{empheq}
\end{subequations}
where $\varGamma(s) = \int_0^{\infty} t^{s-1} e^{-t} dt$ (i.e., the gamma function) and $\gamma := -\varGamma'(1)$ is Euler's constant.

Proposition~\ref{prop_GP} presents the two conditions that are equivalent to Assumption~\ref{assumpt_GEV}.
\begin{prop}\label{prop_GP}[\citealt[Theorems 1.1.3 and 1.1.6]{Haan06}]
Assumption~\ref{assumpt_GEV} holds if and only if each of the conditions (a) and (b) is satisfied:
\begin{description}
\item[\rm{(a)}] There exist some constant $\xi \in \bbR$ and some positive function $a(\,\cdot\,)$ such that
\begin{align}
\lim_{t\to\infty} {U(tx) - U(t) \over a(t)}
=
\dm{x^{\xi} - 1 \over \xi},\qquad x > 0,
\label{connection-a-to-U}
\end{align}
where $(x^{\xi} - 1)/\xi$ is interpreted as $\log x$ if $\xi=0$ (according to the convention introduced after (\ref{defn_G_{xi}})).
\item[\rm{(b)}] For some $\xi \in \bbR$, $F \in \mathsf{MDA}(G_{\xi})$, that is, there exist two sequences of constants $\{ a_n>0;n\in \bb{N}\}$ and $\{ b_n\in\bbR;n\in \bb{N}\}$ such that
\begin{align}
\lim_{n\to\infty}
\PP\!\left( {Z_n - b_n \over a_n} \le z \right) 
= \exp\left\{ -(1+\xi z)_+^{-1/\xi} \right\} 
= G_{\xi}(z),\qquad z \in \bbR.
\label{defn-G_xi}
\end{align}
Note that the constant $\xi$ is the same in (\ref{connection-a-to-U}) and (\ref{defn-G_xi}). In addition, (\ref{defn-G_xi}) holds with $a_n = a(n)$ and $b_n = U(n)$ for $n \in \bbN$.
\end{description}
\end{prop}

\begin{rem}
The standard GEV distribution $G_{\xi}$ is classified into three cases based on the value of $\xi$: (i) $\xi > 0$, (ii) $\xi = 0$, and (iii) $\xi < 0$. The corresponding MDAs for these cases encompass almost all continuous distributions commonly used in statistics (see, e.g., \citealt[Sections 3.3 and 3.4]{Embr97}). Therefore, Assumption~\ref{assumpt_GEV}, which is equivalent to $F \in \mathsf{MDA}(G_{\xi})$, is not overly restrictive.
\end{rem}

\subsection{Regularly varying functions and their scale-freeness in EVT}\label{subsec_class_R}

This subsection presents the results needed for our subsequent discussion on regularly varying functions and their scale-free properties in EVT. First, we introduce the standard definition of regularly varying functions in terms of (asymptotic) scale-freeness. We then summarize the scale-free properties of the functions $a$ and $U$ that appear in Proposition~\ref{prop_GP}; these properties are essential for deriving our power-law formulas for the expected improvement rate and the expected relative gap.

The definition of regularly varying functions is as follows:
\begin{defn}\label{defn_class_R}
An eventually positive and measurable function $f$ on $\bbR_+$ is said to be regularly varying (at infinity) with index $\alpha \in (-\infty,\infty)$ if and only if $f$ is (asymptotically) scale-free, i.e.,
\begin{align}
\lim_{t \to \infty}{f(ct) \over f(t)} = c^{\alpha},
\qquad c > 0.
\label{eqn_decal-free}
\end{align}
The symbol $\calR_{\alpha}$ denotes the set of regularly varying functions with index $\alpha$. In particular, if $f \in \calR_0$, then $f$ is said to be {\it slowly varying}.
\end{defn}

\medskip

\begin{rem}
Strictly speaking, $f$ is scale-free if and only if there exists some $\alpha \in \bbR$ such that $f(ct)=c^{\alpha}f(t)$ for all $t$ in its domain and for all $c>0$. In other words, the basic form of $f$ is invariant under scaling of its argument. Therefore, (\ref{eqn_decal-free}) implies that a regularly varying function is scale-free in an asymptotic sense. For simplicity, we refer to this property of regularly varying functions as ``scale-free."
\end{rem}

\medskip

Proposition~\ref{prop_class_R} presents the fundamental properties of regularly varying functions.
\begin{prop}\label{prop_class_R}
The following statements hold:
\hfill
\begin{enumerate}
\item $f \in \calR_{\alpha}$ if and only if there exists some $L \in \calR_0$ such that $f(t) = L(t) t^{\alpha}$ for $t \in \bbR_+$ (see, e.g., \citealt[Theorem~1.4.1]{Bing89}).
\item If an eventually positive and measurable function $g$ on $\bbR_+$ satisfies $\lim_{t\to\infty} g(t)/f(t) = c$ for some $c > 0$ and $f \in \calR_{\alpha}$, then $g \in \calR_{\alpha}$ (this follows directly from Definition~\ref{defn_class_R}).
\end{enumerate}
\end{prop}

\medskip

\begin{rem}\label{rem_regularly_varying_at_origin}
A positive and measurable function $f$ on $\bbR_+$ is said to be regularly varying at the origin with index $\alpha \in (-\infty,\infty)$, denoted by $\calR_{\alpha}(0)$, if and only if
\begin{align*}
\lim_{t \downarrow 0}{f(ct) \over f(t)} = c^{\alpha},
\qquad c > 0.
\end{align*}
Let $g(t)=f(1/t)$. Then, $g\in \calR_{-\alpha}$ if and only if $f\in \calR_{\alpha}(0)$. Therefore, by Proposition~\ref{prop_class_R}, $f\in \calR_{\alpha}(0)$ if and only if $f(t)=L(1/t)t^{\alpha}$ as $t \downarrow 0$ for some $L\in \calR_0$. This representation shows that $f\in \calR_{\alpha}(0)$ exhibits power-law behavior near the origin.
\end{rem}

\medskip

The asymptotic properties of $a(t)$ and $U(t)$ as $t \to \infty$ play an important role in the analysis of EVT. These properties are summarized in the following proposition (see \citealt[Lemma~1.2.9 and Theorem~B.2.1]{Haan06}).
\begin{prop}\label{prop_function_U}
Suppose that Assumption~\ref{assumpt_GEV} is satisfied, or equivalently (see Proposition~\ref{prop_GP}), that (\ref{connection-a-to-U}) holds for some constant $\xi \in \bbR$ and some positive function $a(\,\cdot\,)$. Under this condition, the following statements are true:
\begin{enumerate}
\item The positive function $a(\,\cdot\,)$ in (\ref{connection-a-to-U}) belongs to the class $\calR_{\xi}$.
\item If $\xi > 0$, then
\begin{align*}
x^* = \lim_{t \to \infty} U(t) &= \infty,
\quad
\lim_{t \to \infty} {U(t) \over a(t) } = {1 \over \xi}.
\end{align*}
\item If $\xi < 0$, then
\begin{align*}
x^* = \lim_{t \to \infty} U(t) &< \infty,
\quad
\lim_{t \to \infty} {x^* - U(t) \over a(t) } = {-1 \over \xi}.
\end{align*}
\item If $\xi = 0$, then $U \in \calR_0$ and
\begin{align}
\lim_{t \to \infty} {a(t) \over U(t)}  = 0.
\label{eqn_240224-01}
\end{align}
Furthermore, if $x^* = \lim_{t \to \infty} U(t) < \infty$, then $U(t) = x^* - L(t)$ for some $L \in \calR_0$ such that $\lim_{t\to\infty}L(t)=0$, and
\begin{align}
\lim_{t \to \infty}  {a(t) \over x^* - U(t)}&= 0.
\label{eqn_240224-02}
\end{align}
\end{enumerate}
 
\end{prop}

\section{Power-law formulas behind the curse of scale-freeness}\label{sec_main}

This section discusses the power-law behavior of empirical solutions generated by the RMS method. First, we present the main results of this paper, the {\it power-law formulas}, which establish the power-law properties of the expected improvement rate $\EE[R_n(x)]$ and the expected relative gap $\EE[\varDelta_n(x)]$. We then derive the simplified versions of these formulas to highlight the power-law characteristics in $\EE[R_n(x)]$ and $\EE[\varDelta_n(x)]$. Finally, we validate these theoretical results through numerical experiments (see Section~\ref{sec_num_setup} for the experimental setup).

Theorem~\ref{thm_asymp-E[R_n(x)]-01} presents the power-law formulas for $\EE[R_n(x)]$ and $\EE[\varDelta_n(x)]$.
\begin{thm}\label{thm_asymp-E[R_n(x)]-01}
Suppose that $\EE[X] =\int_0^{\infty} x \,d F(x) < \infty$ and Assumption~\ref{assumpt_GEV} holds for $G=G_{\xi}$ with $\xi < 1$. Then, the following statements are true:
\begin{enumerate}
\item If $0 < \xi < 1$,
\begin{align}
\lim_{n\to\infty} {1 \over a(n)}\EE[R_n(x)]
&= -{\varGamma(-\xi) \over x},
 \qquad x > 0,
\label{lim_E[R_n(x)]_xi>0}
\end{align}
where $a \in \calR_{\xi}$ (due to Proposition~\ref{prop_function_U}(i)).
\item If $\xi = 0$ and $x^* = \infty$, then
\begin{align}
\lim_{n\to\infty} {1 \over U(n)}\EE[R_n(x)]
&= {1 \over x},
\qquad x > 0,
\label{lim_E[R_n(x)]_xi=0_infinite}
\end{align}
where $U \in \calR_0$ (due to Proposition~\ref{prop_function_U}(iv)).
\item If $\xi = 0$ and $x^* < \infty$, then
\begin{subequations}\label{power-law_(iii)}
\begin{alignat}{2}
\lim_{n\to\infty} {1 \over x^* - U(n)}
\left({x^* - x \over x} -  \EE[R_n(x)]\right)
&= {1 \over x},
&\qquad 0 &< x < x^*,
\label{lim_E[R_n(x)]_xi=0_finite}
\\
\lim_{n\to\infty} {1 \over x^* - U(n)} \EE[\varDelta_n(x)]
&= {1 \over x^* },
&\qquad 0 &< x < x^*,
\label{lim_E[Delta_n(x)]_xi=0_finite}
\end{alignat}
\end{subequations}
where $U(t) = x^* - L(t)$ for some $L \in \calR_0$ such that $\lim_{t\to\infty}L(t)=0$ (due to Proposition~\ref{prop_function_U}(iv)).
\item If $\xi < 0$, then
\begin{subequations}\label{power-law_(iv)}
\begin{alignat}{2}
\lim_{n\to\infty} {1 \over a(n)}\left({x^* - x \over x} - \EE[R_n(x)]\right)
&= {\varGamma(-\xi) \over x},
&\qquad 0 &< x < x^*,
\label{lim_E[R_n(x)]_xi<0}
\\
\lim_{n\to\infty} {1 \over a(n)} \EE[\varDelta_n(x)]
&= {\varGamma(-\xi) \over x^* },
&\qquad 0 &< x < x^*,
\label{lim_E[Delta_n(x)]_xi<0}
\end{alignat}
\end{subequations}
where $a \in \calR_{\xi}$ (due to Proposition~\ref{prop_function_U}(i)).
\end{enumerate}

\end{thm}

\begin{proof}
See Appendix~\ref{proof_thm_asymp-E[R_n(x)]-01}. 
\end{proof}

\begin{rem}
Equations (\ref{power-law_(iii)}) and (\ref{power-law_(iv)}) reflect the affine relation, expressed in (\ref{eqn_affine_relation}), between $R_n(x)$ and $\varDelta_n(x)$.
\end{rem}

\medskip

\begin{rem}
We confirm that the right-hand sides of  (\ref{lim_E[R_n(x)]_xi>0}), (\ref{lim_E[R_n(x)]_xi<0}), and (\ref{lim_E[Delta_n(x)]_xi<0}) are positive because $\varGamma(-\xi) > 0$ for $\xi < 0$ and $-\varGamma(-\xi)>0$ for $0 < \xi < 1$. The former follows from the definition of the gamma function $\varGamma$, and the latter follows from Euler's reflection formula, which states that $\varGamma(z) \varGamma(1 - z) = \pi / \sin (\pi z)$ for $z \not\in \{0,\pm1,\pm2,\dots\}$.
\end{rem}

\medskip

Theorem~\ref{thm_asymp-E[R_n(x)]-01}, in conjunction with Proposition~\ref{prop_class_R}(i), yields the simplified versions of the power-law formulas for $\EE[R_n(x)]$ and $\EE[\varDelta_n(x)]$. This simplification clarifies the power-law characteristics of the two expectations as functions of the number of trials.
\begin{coro}\label{coro_regularly-varying}
Under the conditions of Theorem~\ref{thm_asymp-E[R_n(x)]-01}, the following statements are true:
\begin{enumerate}
\item If $0 \le \xi < 1$ and $x^* = \infty$, then
\begin{align*}
\EE[R_n(x)] = L(n)n^{\xi},\qquad  x > 0,
%\label{ean_EIR_in_R_{xi}_infinite}
\end{align*}
for some $L \in \calR_0$.
\item If $\xi \le 0$ and $x^* < \infty$, then
\begin{alignat}{2}
\EE[R_n(x)] 
&= {x^* - x \over x} - {x^*  \over x} L(n)n^{\xi},
&\qquad 0 &< x < x^*,
\nonumber
%\label{ean_EIR_in_R_{xi}_finite}
\\
\EE[\varDelta_n(x)] &= L(n)n^{\xi},
&\qquad 0 &< x < x^*,
\label{ean_ERG_in_R_{xi}}
\end{alignat}
for some $L \in \calR_0$.
\end{enumerate}
\end{coro}

The simplified formulas in Corollary~\ref{coro_regularly-varying} yield Proposition~\ref{prop_scale-free-01} below, which supports the perspective introduced in Section~\ref{sec_introduction}: the evolution of the relative gap of the best EOV is subject to the {\it curse of scale-freeness} (see (\ref{eqn_curse_of_scale_freenss})) and is approximately linear on a log-log scale as a function of the number of trials $n$ (see Figure~\ref{fig_RE-NN-LKH-log}).
\begin{prop}\label{prop_scale-free-01}
Suppose that the conditions of Theorem~\ref{thm_asymp-E[R_n(x)]-01} hold. Furthermore, suppose that $\xi \le 0$ and $x^* < \infty$. We then have
\begin{alignat}{2}
\lim_{n\to\infty} {\EE[\varDelta_{cn}(x)] \over \EE[\varDelta_n(x)]}
&=  c^{\xi}, &\qquad 0 &< x < x^*,~c \in \bbN,
\label{scale-free-ERG-01}
\\
\lim_{n\to\infty}{ \log \EE[\varDelta_n(x)] \over \log n} 
&= \xi, &\qquad  0 &< x < x^*.
\label{scale-free-ERG-02}
\end{alignat}
\end{prop}

\begin{proof}
See Appendix~\ref{proof:prop_scale-free-01}
\end{proof}

Proposition~\ref{prop_scale-free-01} implies the curse of scale-freeness, which characterizes the intractability of large-scale optimization with multi-start methods. As mentioned in Section~\ref{sec_introduction}, Equation~(\ref{scale-free-ERG-01}) leads to (\ref{eqn_curse_of_scale_freenss}), indicating that the number of additional trials required to halve the expected relative gap is approximately $(2^{-1/\xi} - 1)$ times the number of completed trials. This approximation is likely to become more accurate as the number of completed trials increases. It would be reasonable to describe such a situation figuratively as ``{\it Reaching for the goal makes it slip away}."

\begin{figure*}[htbp]
\begin{tabular}{cc}
\begin{minipage}[t]{0.48\hsize}
\centering
\includegraphics[scale=0.45,bb=0 0 444 316]{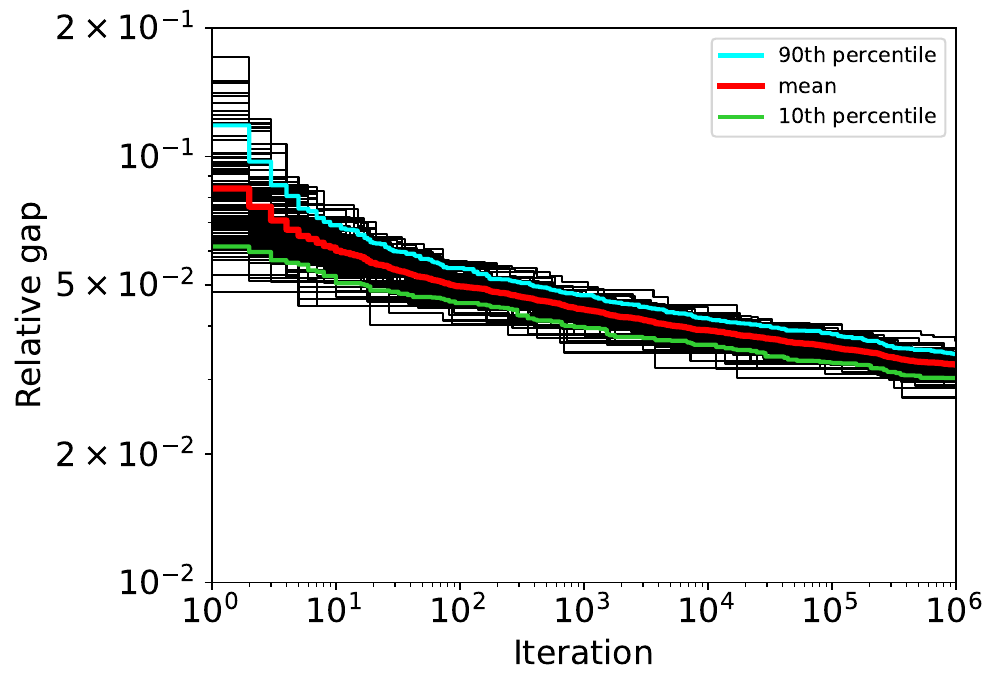}
\subcaption{RA + 3-opt}
\label{fig_rel_gap-RA-3opt}
\end{minipage} &
\begin{minipage}[t]{0.48\hsize}
\centering
\includegraphics[scale=0.45,bb=0 0 444 316]{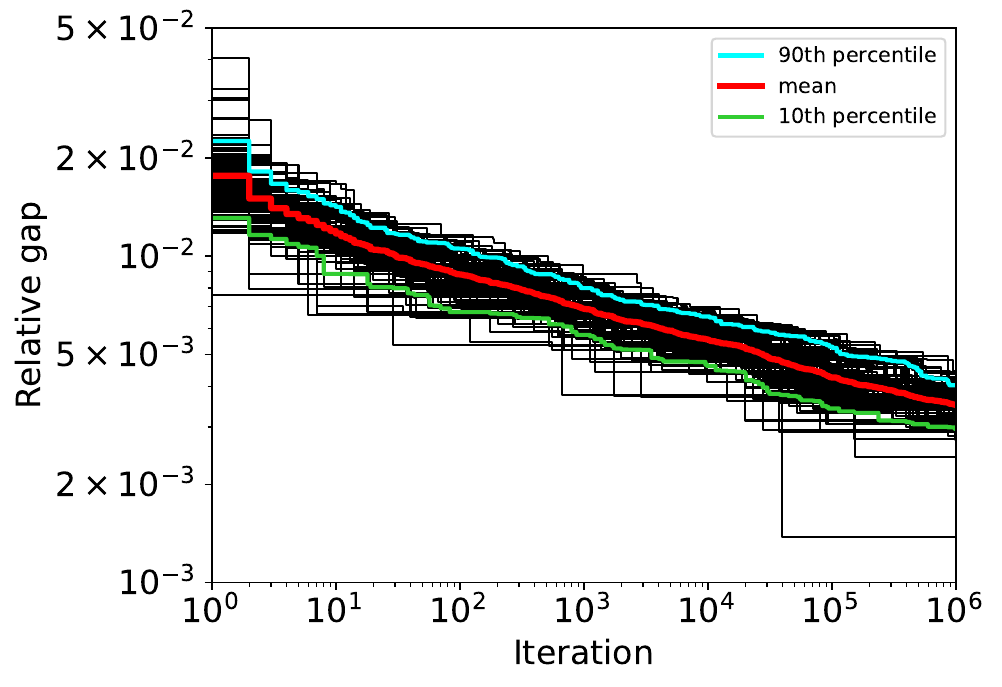}
\subcaption{RA + LK}
\label{fig_rel_gap-RA-LKH}
\end{minipage} \\[10mm]
\begin{minipage}[t]{0.48\hsize}
\centering
\includegraphics[scale=0.45,bb=0 0 444 316]{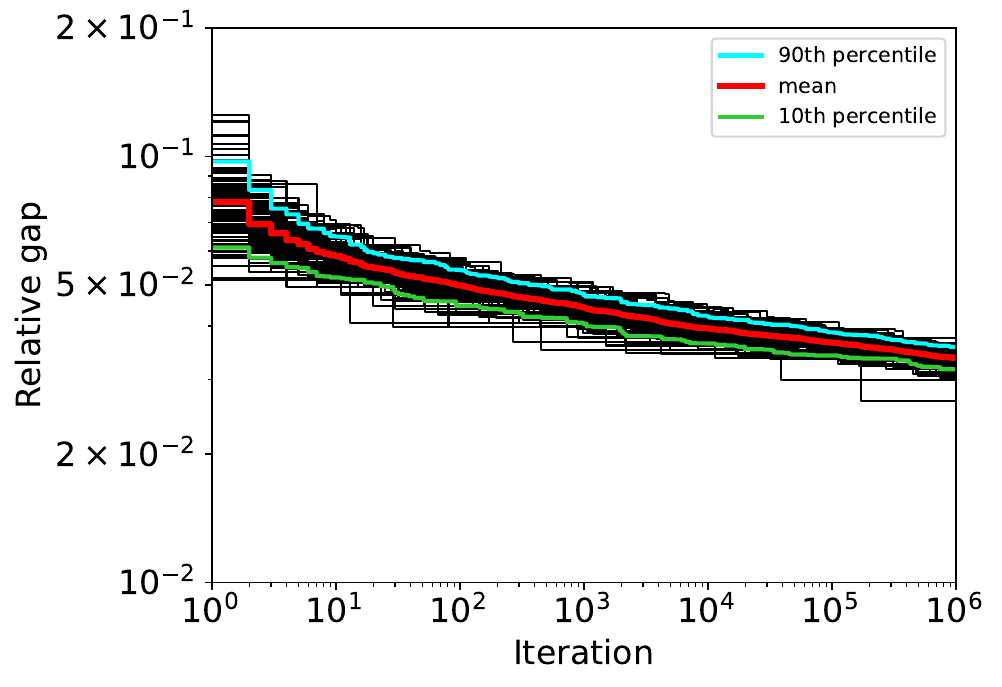}
\subcaption{NN + 3-opt}
\label{fig_rel_gap-NN-3opt}
\end{minipage} &
\begin{minipage}[t]{0.48\hsize}
\centering
\includegraphics[scale=0.45,bb=0 0 444 316]{fig2_RMS_3/NN_LK_p001-100_try001_xlog_ylog_rel_gap_to_true_opt.pdf}
\subcaption{NN + LK}
\label{fig_rel_gap-NN-LKH}
\end{minipage} \\[10mm]
\begin{minipage}[t]{0.48\hsize}
\centering
\includegraphics[scale=0.45,bb=0 0 444 316]{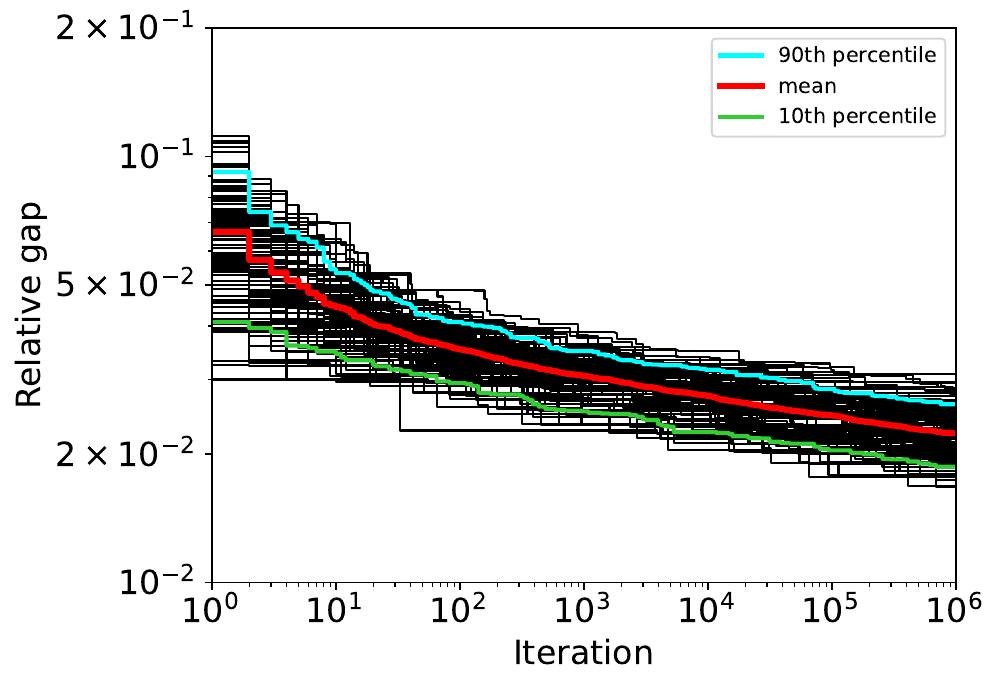}
\subcaption{GR + 3-opt}
\label{fig_rel_gap-GR-3opt}
\end{minipage} &
\begin{minipage}[t]{0.48\hsize}
\centering
\includegraphics[scale=0.45,bb=0 0 444 316]{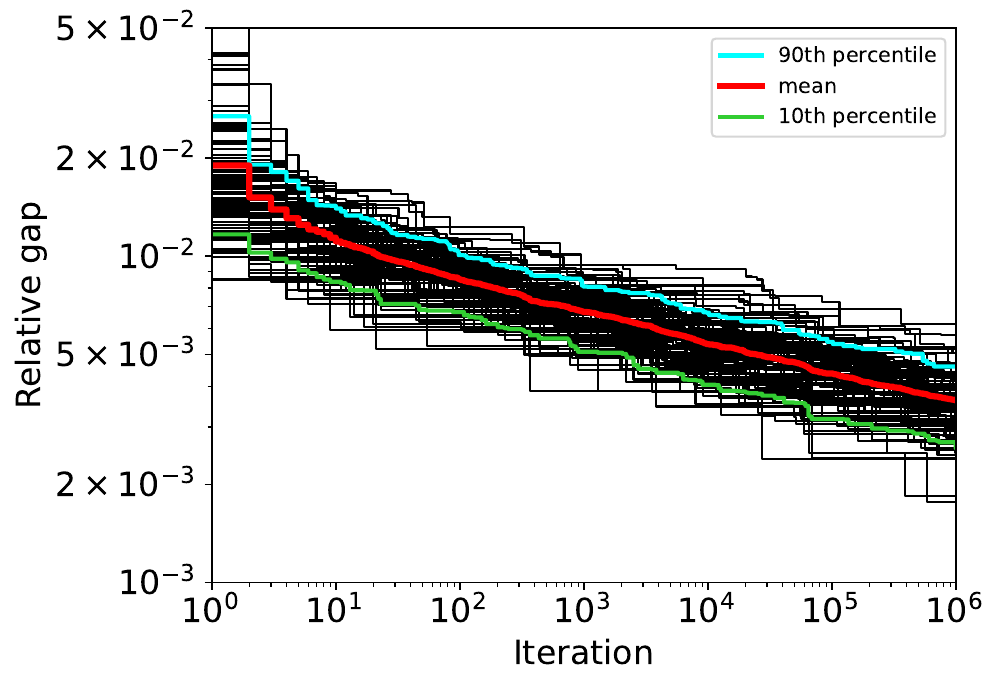}
\subcaption{GR + LK}
\label{fig_rel_gap-GR-LKH}
\end{minipage}
\end{tabular}
\caption{Evolution of the relative gap of the best EOV over $10^6$ iterations for each of the six RMS algorithms listed in Table~\ref{tb_RMS_algorithms}, applied to 100 random TSP instances. The blue, red, and green lines represent the 90th percentile, mean, and 10th percentile values, respectively, based on the results obtained from the 100 instances.}\label{fig_RMS_1000cities}
\end{figure*}

Finally, we perform numerical experiments to confirm the occurrence of the curse of scale-freeness in large-scale optimization with multi-start methods. To achieve this, we solve 100 random instances of the TSP with 1000 cities using six different RMS algorithms (see Table~\ref{tb_RMS_algorithms}). For each instance, each algorithm generates $10^6$ locally optimal solutions (some of which are duplicates). Figure~\ref{fig_RMS_1000cities} shows the log-log plots of the evolution of the relative gap of the best EOV as the six RMS algorithms are applied to 100 random TSP instances (see Appendix~\ref{subsec_how_to_run_algorithms} for details on running the algorithms to generate this figure). The results imply that the relative gap of the best EOV decreases with the number of iterations at a power-law rate, supporting the presence of the curse of scale-freeness.

%%%%%%%%%%%%%%%%%%%%%%%%%%%%%% Figure 3 %%%%%%%%%%%%%%%%%%%%%%%%%%%%%%

\section{Understanding the curse of scale-freenes}\label{sec_discussion}

This section discusses the curse of scale-freeness in large-scale optimization from two perspectives. Section~\ref{subsec_good_solutions} analyzes how this curse is linked to the ratio of good solutions, while Section~\ref{subsec_intractability} addresses how challenging it is to overcome the curse.

\subsection{The curse of scale-freeness in connection to the ratio of good solutions}\label{subsec_good_solutions}

This subsection discusses the connection between the curse of scale-freeness and the ratio of good (empirical) solutions. First, we introduce the good-solution-ratio function $r(\varep)$ for $\varep \in (0,1)$. We then prove that the curse of scale-freeness arises from the scale-freeness of the left tail of $r(\varep)$, that is, its power-law-type decay as $\varep$ approaches zero. This theoretical result is supported by numerical experiments. Furthermore, we show that the situation worsens if the left tail of $r(\varep)$ decays faster than any power-law rate.

We introduce a function related to the ratio of good solutions generated by the RMS method. Assuming that the supremum $x^*$ of the EOVs is finite, we define $r(\varepsilon)$ for $0 < \varepsilon < 1$ as follows:
\begin{align}
r(\varep) = \PP \left( \frac{x^* - X}{x^* } < \varep \right)
= \PP(X > x^*- \varep x^* ) = 1 - F(x^*- \varep x^*),
\label{defn_r(varep)}
\end{align} 
where $\varep \in (0,1)$ is generally assumed to be small. Since the EOVs $X_n$, $n \in \bbN$, are i.i.d., $r(\varep)$ can be interpreted as the ratio of good solutions whose objective values are within a relative gap (error) $\varep$ from the supremum $x^*$ of the EOVs. Thus, $r(\varep)$ is referred to as the {\it good-solution-ratio function}. With the good-solution-ratio function $r$, we can express the distribution $F$ of EOVs as follows:
\begin{align}
F(x) = 1 - r\left({x^* - x \over x^*} \right),\qquad 0 \le x < x^*.
\label{eqn:F_expressed_by_r}
\end{align}

The following proposition relates the scale-freeness of the good-solution-ratio function $r(\varep)$ (as $\varep \downarrow 0$) to that of the expected relative gap $\EE[\varDelta_n(x)]$ (as $n \to \infty$).
\begin{prop}\label{prop_power-law-decay} 
Suppose that $x^* < \infty$ and thus $\EE[X] =\int_0^{x^*} x \,d F(x) < \infty$. Furthermore, if there exist some $\psi > 0$ and $L \in \calR_0$ such that $r(\varep) = L(1/\varep)\varep^{\psi}$, or equivalently,
\begin{align}
\lim_{\varep \downarrow 0}{r(c\varep) \over r(\varep)} = c^{\psi},
\qquad c > 0,
\label{scale-Free-r(varep)}
\end{align}
then $F \in \mathsf{MDA}(G_{\xi})$ with $\xi = -1/\psi < 0$, and thus both (\ref{scale-free-ERG-01}) and (\ref{scale-free-ERG-02}) hold.
\end{prop}

\begin{proof}
See Appendix~\ref{proof:prop_power-law-decay}
\end{proof}

\begin{rem}
The scale-freeness property, expressed in (\ref{scale-Free-r(varep)}), of the good-solution-ratio function $r$ implies that the distribution $F$ of EOVs exhibits a power-law behavior near its right endpoint $x^*$. Indeed, combining (\ref{eqn:F_expressed_by_r}) with $r(\varep) = L(1/\varep)\varep^{\psi}$ yields
\begin{align*}
1 - F(x) = L\left( {x^* \over x^* - x} \right)
\left( {x^* - x \over x^*} \right)^{\psi},
\qquad \mbox{as $x \uparrow x^*$}.
\end{align*}
\end{rem}

Proposition~\ref{prop_power-law-decay} provides a valuable insight into the curse of scale-freeness: regardless of how small the current relative gap $\varep$ (from the supremum $x^*$) of the best empirical solution becomes, the fraction of solutions within a relative gap of $\varep/2$ asymptotically remains constant. In fact, (\ref{scale-Free-r(varep)}) yields
\begin{align*}
\lim_{\varep \downarrow 0}{r(\varep/2) \over r(\varep)} = \left( {1 \over 2} \right)^{\psi}.
\end{align*}
This means that, asymptotically, the probability of finding a solution within a relative gap of $\varep/2$ remains a constant fraction of the probability of finding one within a relative gap of $\varep$. In short, the curse of scale-freeness implies that---\textit{Even if you find a good solution, many better ones may still exist.} To support this insight, Figure~\ref{fig_GSR_RMS_1000cities}, based on the same data as Figure~\ref{fig_RMS_1000cities}, plots the good-solution-ratio function $r(\varep)$, $\varep \in (0,1)$, on a log-log scale and shows that it becomes linear as $\varep$ decreases toward zero, indicating power-law behavior. These results support Proposition~\ref{prop_power-law-decay}.

To complement Proposition~\ref{prop_power-law-decay}, we present another proposition that addresses the case where the good-solution-ratio function  $r(\varep)$ has a faster left-tail decay rate---specifically, an exponential rate (e.g., $\exp\{-\varep^{-\alpha}\}$ with $\alpha > 0$).
\begin{prop}\label{prop_exponential-type-decay} 
Suppose that $x^* < \infty$ and the good-solution-ratio function $r$ has the following representation: 
\begin{align}
r(\varep) = c_0\exp\{-g(1/\varep)\},\qquad 0 < \varep < 1,
\label{cond_F_moderate_exp}
\end{align}
where $c_0>0$ is a constant and $g: \bbR_+ \to (0,\infty)$ is a twice differentiable function such that $g \in \calR_{\alpha}$ for some $\alpha > 0$ and
\begin{align}
\lim_{t \to \infty} tg'(t) &= \infty,
\label{cond_g(t)-01}
\\
\lim_{t \to \infty} {g''(t) \over [g'(t)]^2 } &= 0.
\label{cond_g(t)-02}
\end{align}
We then have $F \in \mathsf{MDA}(G_0)$ and
\begin{align}
\EE[\varDelta_n(x)]
&= {x^* - U(n) \over x^* } + o(x^* - U(n)),
\qquad \mbox{as $n \to \infty$},
\label{asymp_Delta_n(x)}
\end{align}
for $x \in (0, x^*)$, where $U(t) = x^* - L(t)$ for some $L \in \calR_0$ such that $\lim_{t\to\infty}L(t)=0$.
\end{prop}

\begin{proof}
See Appendix~\ref{proof:prop_exponential-type-decay}. 
\end{proof}

Proposition~\ref{prop_exponential-type-decay} suggests that the difficulty of closing the expected relative gap $\EE[\varDelta_n(x)]$ of the best EOV increases when the good-solution-ratio function $r$ has a faster left-tail decay. This increased difficulty is because, as (\ref{asymp_Delta_n(x)}) shows, $\EE[\varDelta_n(x)]$ is slowly varying as a function of $n$, and thus its decay is slower than any power function. A typical example of this case is as follows:
\begin{align}
\lim_{n\to\infty} {\EE[\varDelta_n(x)] \over (\log n)^{-\beta}} = c_0
\quad \mbox{for some $\beta > 0$ and $c_0 > 0$}.
\label{eqn:add_250321-01}
\end{align}
Such behavior of $\EE[\varDelta_n(x)]$ is even more problematic than the curse of scale-freeness. This difficulty likely arises from the challenge of finding {\it better} solutions among a sparse set of good solutions scattered across a large feasible domain.

To further illustrate the above discussion, we consider a specific example from the cases covered by Proposition~\ref{prop_exponential-type-decay} and derive an  explicit limit formula for the expected relative gap $\EE[\varDelta_n(x)]$. Assume that
\begin{align}
r(\varep) = c_0 \exp\{-\phi \varep^{-\alpha}\},\qquad
0 < \varep < 1,
\label{examp_r(verep)}
\end{align}
where $c_0 > 0$, $\phi > 0$, and $\alpha > 0$. Substituting (\ref{examp_r(verep)}) into (\ref{eqn:F_expressed_by_r}), we have
\begin{align*}
1 - F(x) 
= c_0 \exp \left\{-\phi\left( {x^* \over x^* - x} \right)^{\alpha} \right\},
\qquad 0 \le x < x^*.
\end{align*}
From this equation and (\ref{defn_U}), we obtain
\begin{align}
{x^* - U(t) \over x^*}
= [\phi^{-1} \log (c_0t) ]^{-1/\alpha}.
\label{eqn_U(x)}
\end{align}
Applying (\ref{eqn_U(x)}) to (\ref{asymp_Delta_n(x)}) yields
\begin{align}
\EE[\varDelta_n(x)]
&= [\phi^{-1} \log (c_0 n) ]^{-1/\alpha}
+ o((\log n)^{-1/\alpha}),
\qquad \mbox{as $n \to \infty$},
\label{asymp_examp_Delta_n(x)}
\end{align}
for $x \in (0,x^*)$, and thus (\ref{eqn:add_250321-01}) holds with $\beta = 1/\alpha$. Furthermore, (\ref{asymp_examp_Delta_n(x)}) shows that the decay of the expected relative gap $\EE[\varDelta_n(x)]$ is dominated by $(\log n)^{-1/\alpha}$, resulting in an extremely slow rate. As mentioned above, in this scenario, the left tail of the good-solution-ratio function $r$ decays even slower than any power-law function.

Building on the above analysis, we conclude that the {\it curse of scale-freeness}---expressed by the metaphor ``{\it Reaching for the goal makes it slip away}''---inevitably arises in large-scale optimization when using RMS methods, regardless of the decay rate of the good-solution-ratio function $r$. Notably, if the left tail of $r$ decays faster than any power-law rate, the resulting intractability can be more severe than in typical cases of the curse of scale-freeness.

\subsection{Intractability of overcoming the curse of scale-freeness}\label{subsec_intractability}

This subsection explores the curse of scale-freeness---that is, the intractability of large-scale optimization with multi-start methods. We first show that overcoming this curse requires a powerful LS algorithm using effective restart and diversification strategies to exponentially accelerate solution improvement relative to the RMS method. We then present numerical experiments demonstrating that even the ILS method \citep{Lou19}, a widely recognized and powerful metaheuristic approach, is unlikely to overcome the curse of scale-freeness. Finally, we summarize our theoretical and numerical results, highlighting what the curse of scale-freeness is and how it arises in large-scale optimization with multi-start methods.

The power-law formula (\ref{ean_ERG_in_R_{xi}}) illustrates how difficult it is to design a metaheuristic algorithm capable of overcoming the curse of scale-freeness. To illustrate this difficulty, we consider two types of acceleration for the RMS method: one is \emph{polynomial} acceleration and the other is \emph{exponential} acceleration. The expected relative gaps for polynomially and exponentially accelerated RMS methods are obtained by replacing $n$ in (\ref{ean_ERG_in_R_{xi}}) with $n^K$ and $K^n$ ($K = 2, 3, \dots$), respectively:
\begin{subequations}\label{two_accelerations}
\begin{alignat}{3}
&\mbox{(Polynomial acceleration)} &\quad
\EE[\varDelta_{n^K}(x)] 
&= L(n^K) n^{K\xi},
& \qquad &\mbox{$\xi < 0$},
\label{Polynomial_case}
\\
&\mbox{(Exponential acceleration)} &\quad
\EE[\varDelta_{K^n}(x)] 
&=  L(K^n) K^{n\xi},
& \qquad &\mbox{$\xi < 0$},
\label{Exponential_case}
\end{alignat}
\end{subequations}
where $L$ is a slowly varying function that may differ from the $L$ in (\ref{ean_ERG_in_R_{xi}}). Equation (\ref{two_accelerations}) shows that the scale-free property persists under polynomial acceleration but not under exponential acceleration. Thus, overcoming the curse of scale-freeness requires a powerful LS algorithm---one that, through effective restart and diversification strategies, exponentially accelerates solution improvement relative to the {\it standard} RMS method. However, achieving such exponential acceleration while keeping the cost per trial manageable is challenging without a sufficiently tractable problem structure.

\begin{figure*}[htbp]
 \begin{tabular}{l@{~~}cc}
  \hline
  & RMS & ILS \\
  \hline
  \verb|brd14051|
  &
  \begin{minipage}{0.42\hsize}
   \centering
   \includegraphics[width=5.5cm, bb=0 0 444 317]{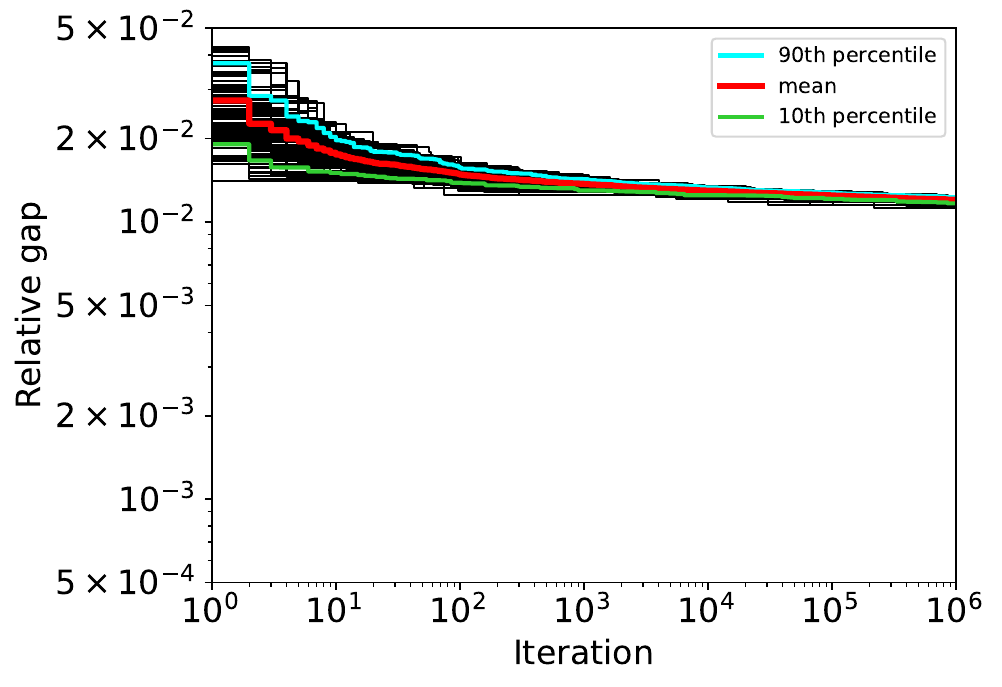}
  \end{minipage}    
  &
  \begin{minipage}{0.42\hsize}
   \centering
   \includegraphics[width=5.5cm, bb=0 0 444 317]{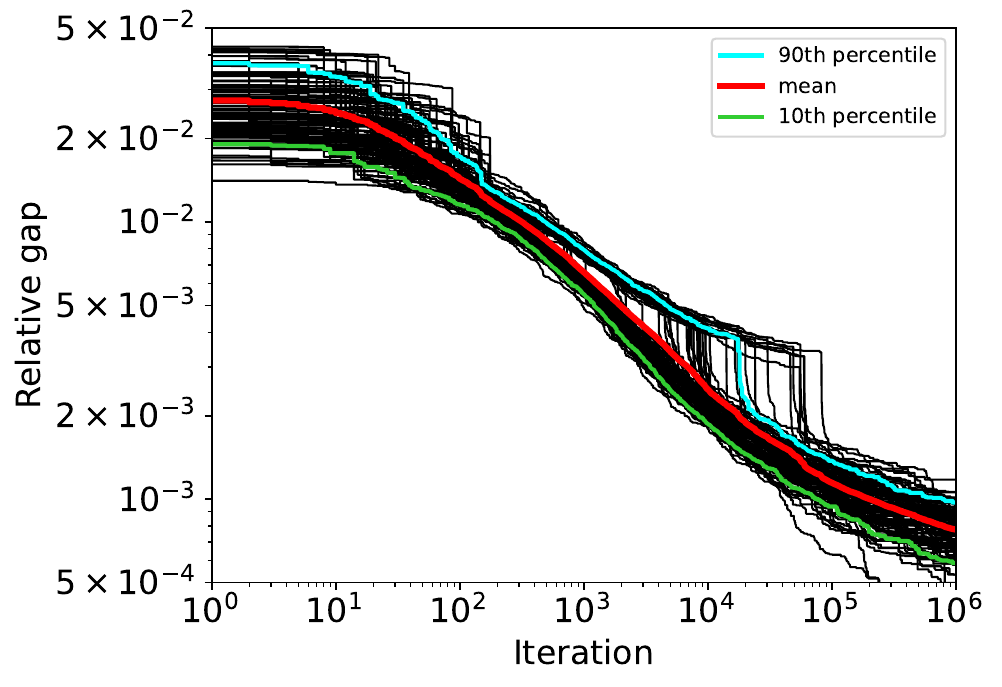}
  \end{minipage}    
  \\
  \hline
  \verb|d15112|
  &
  \begin{minipage}{0.42\hsize}
   \centering
   \includegraphics[width=5.5cm, bb=0 0 444 317]{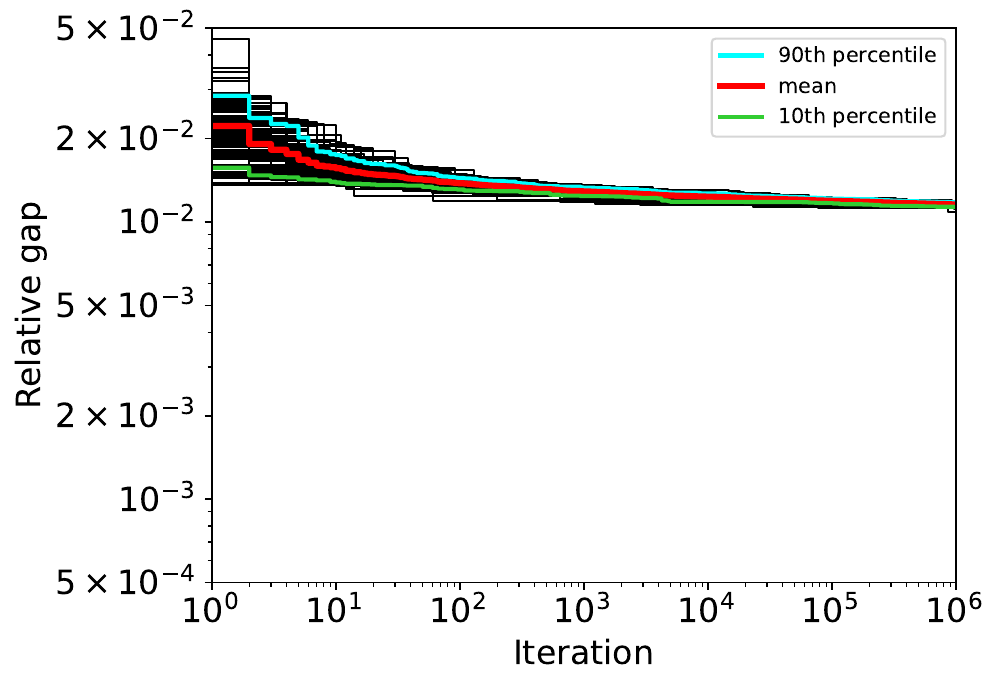}
  \end{minipage}    
  &
  \begin{minipage}{0.42\hsize}
   \centering
   \includegraphics[width=5.5cm, bb=0 0 444 317]{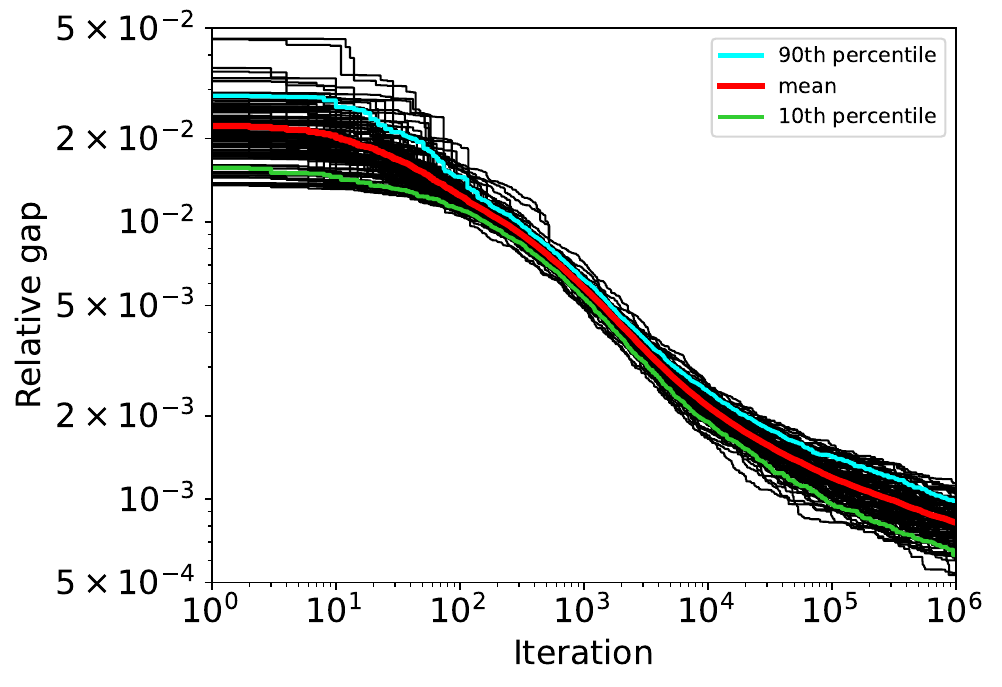}
  \end{minipage}    
  \\
  \hline
  \verb|d18512|
  &
  \begin{minipage}{0.42\hsize}
   \centering
   \includegraphics[width=5.5cm, bb=0 0 444 317]{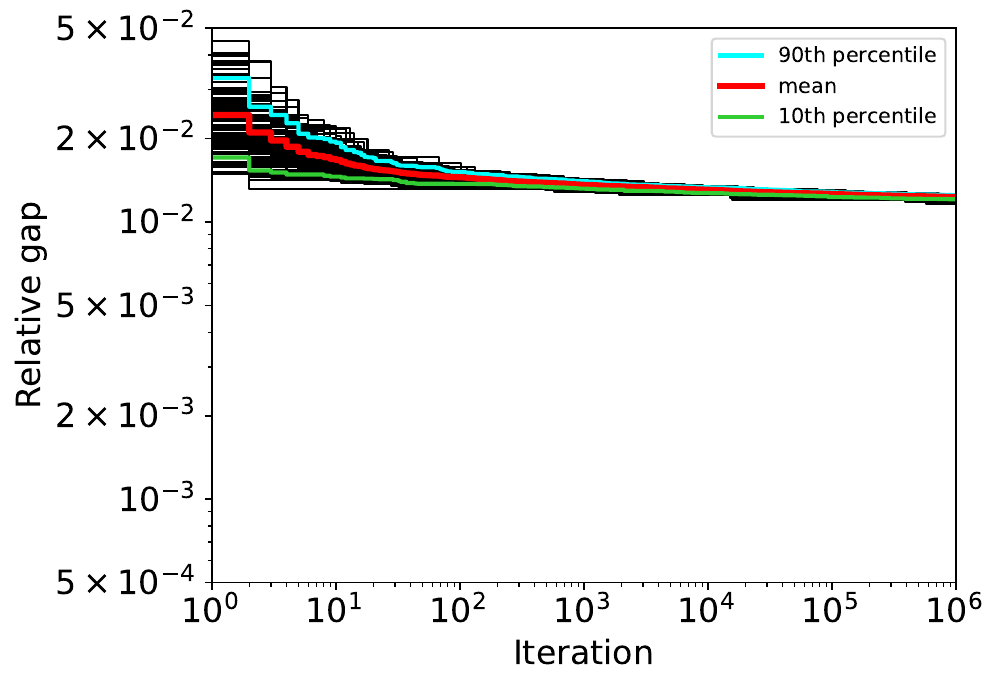}
  \end{minipage}    
  &
  \begin{minipage}{0.42\hsize}
   \centering
   \includegraphics[width=5.5cm, bb=0 0 444 317]{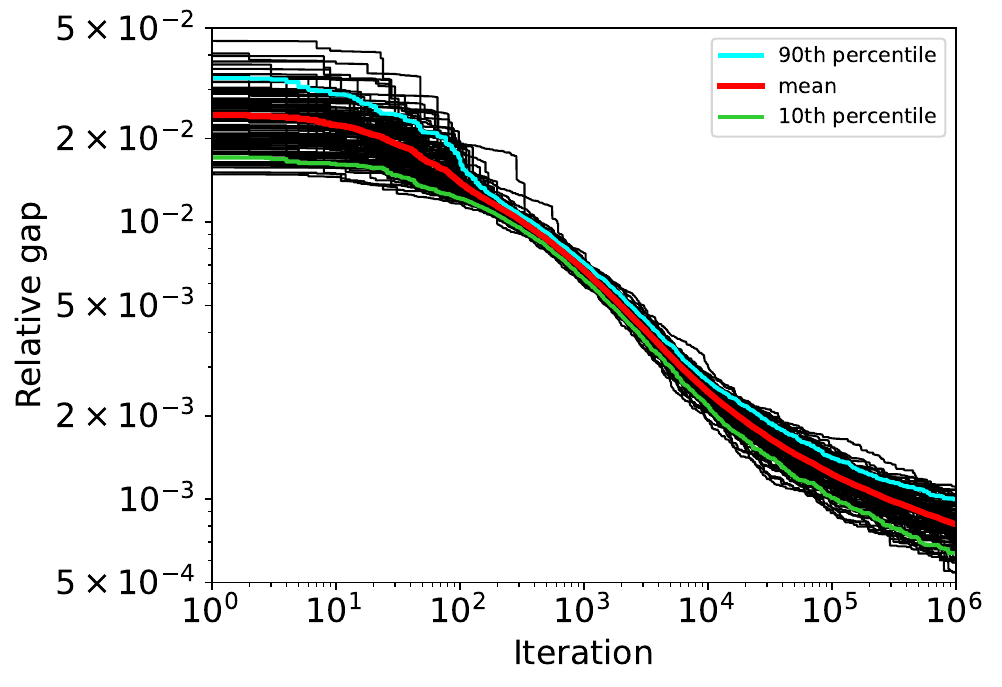}
  \end{minipage}    
  \\
  \hline
  \verb|rl11849|
  &
  \begin{minipage}{0.42\hsize}
   \centering
   \includegraphics[width=5.5cm, bb=0 0 444 317]{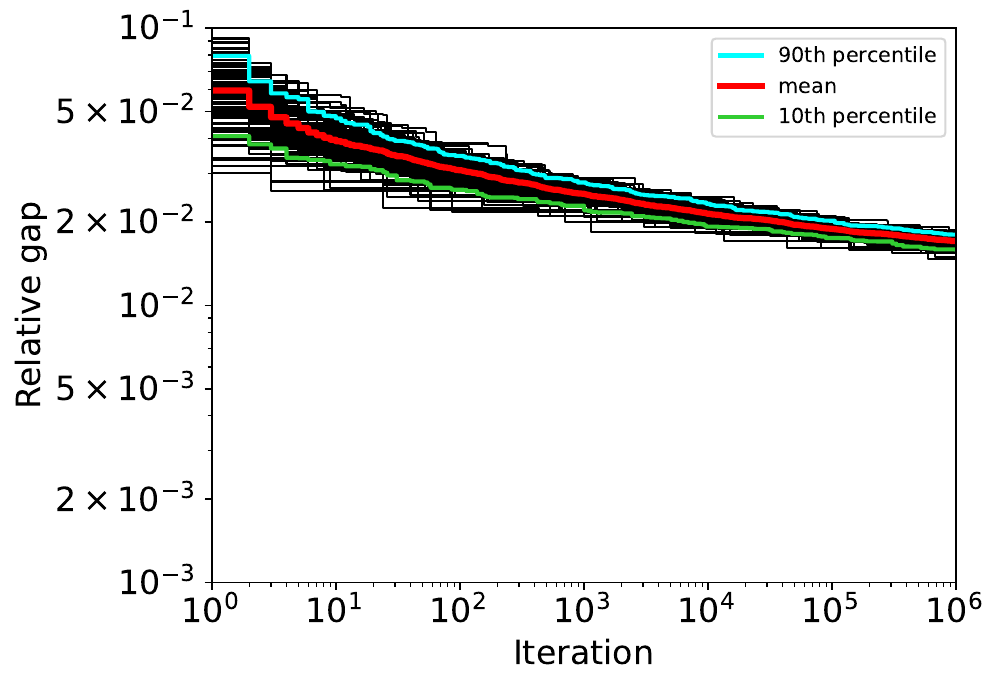}
  \end{minipage}    
  &
  \begin{minipage}{0.42\hsize}
   \centering
   \includegraphics[width=5.5cm, bb=0 0 444 317]{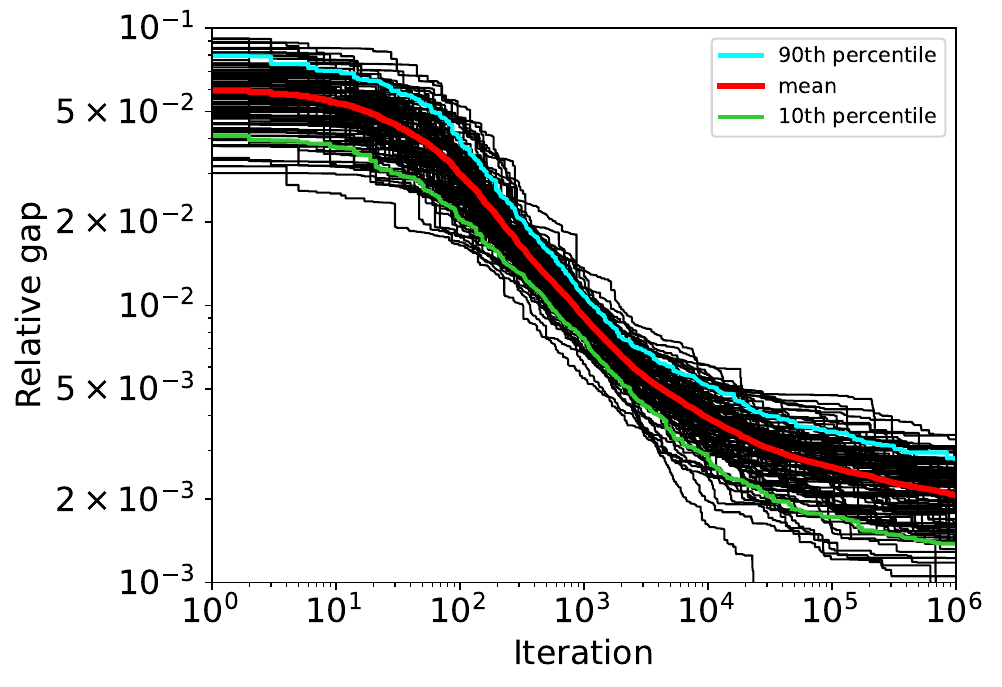}
  \end{minipage}    
  \\
  \hline
  \verb|usa13509|
  &
  \begin{minipage}{0.42\hsize}
   \centering
   \includegraphics[width=5.5cm, bb=0 0 444 317]{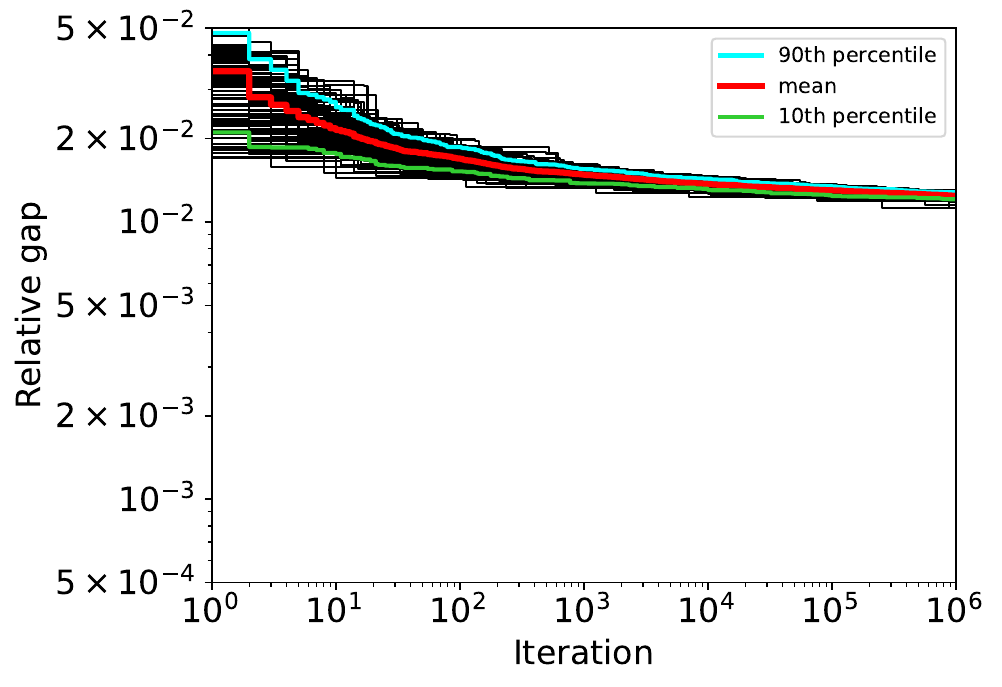}
  \end{minipage}    
  &
  \begin{minipage}{0.42\hsize}
   \centering
   \includegraphics[width=5.5cm, bb=0 0 444 317]{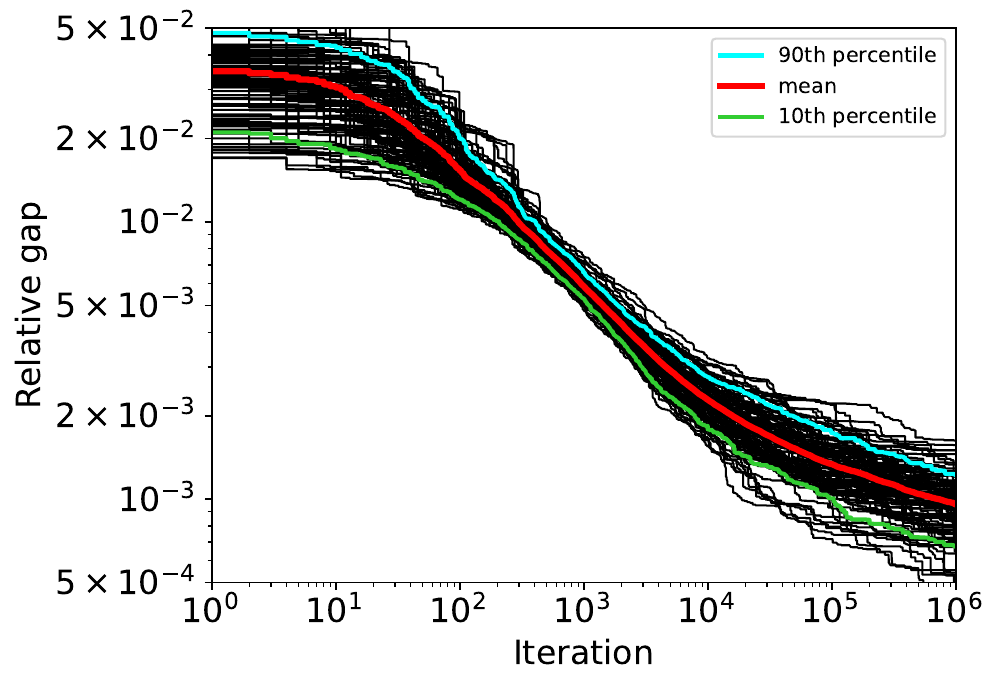}
  \end{minipage}    
  \\
  \hline
 \end{tabular}
 \caption{Evolution of the relative gap of the best EOV over 100 random sets of $10^6$ iterations of the RA+LK RMS and ILS algorithms, respectively, for the five TSPLIB instances listed in Table~\ref{TSPLIB_instances}.}\label{RMS-ILS_RA-LKH}
\end{figure*}

\begin{figure*}[htbp]
 \begin{tabular}{ccc}
  \hline
  & RMS & ILS \\
  \hline
  \verb|brd14051|
  &
  \begin{minipage}{0.42\hsize}
   \centering
   \includegraphics[width=5.5cm, bb=0 0 444 317]{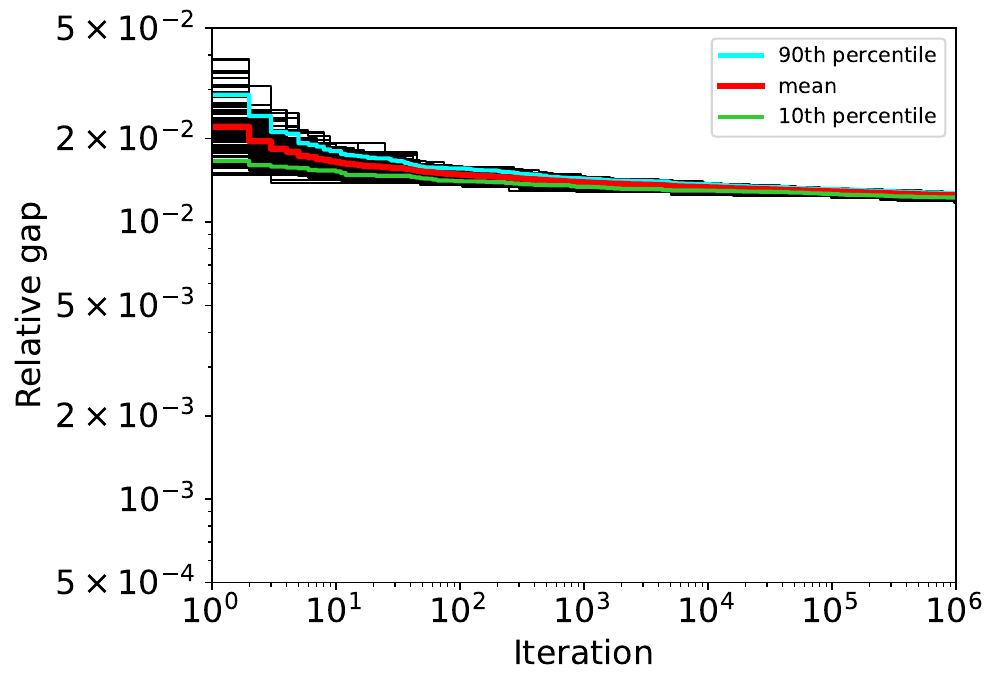}
  \end{minipage}    
  &
  \begin{minipage}{0.42\hsize}
   \centering
   \includegraphics[width=5.5cm, bb=0 0 444 317]{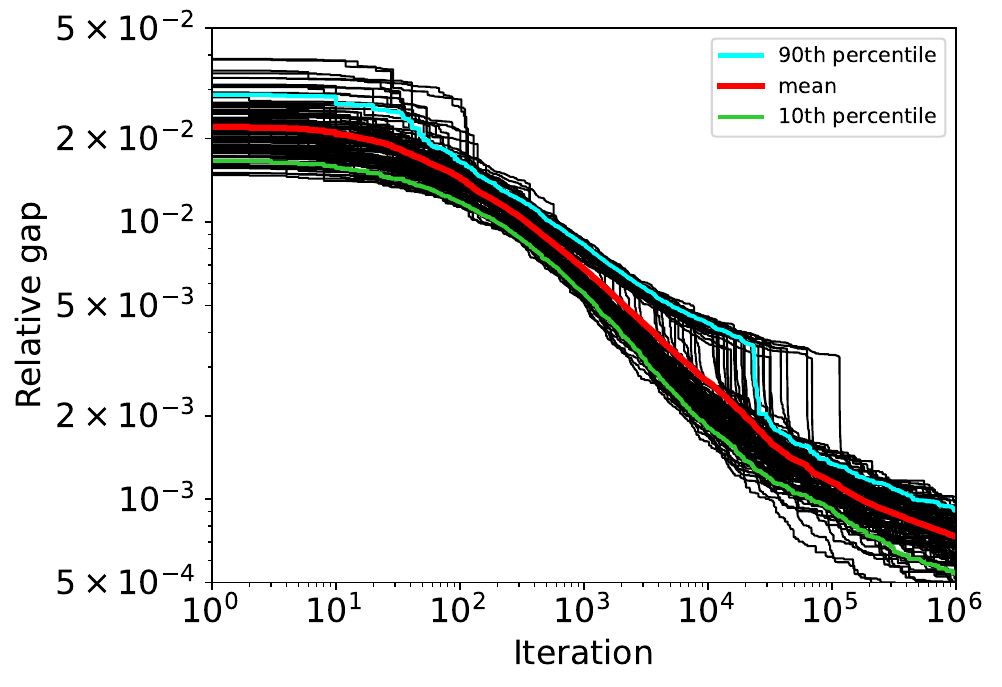}
  \end{minipage}    
  \\
  \hline
  \verb|d15112|
  &
  \begin{minipage}{0.42\hsize}
   \centering
   \includegraphics[width=5.5cm, bb=0 0 444 317]{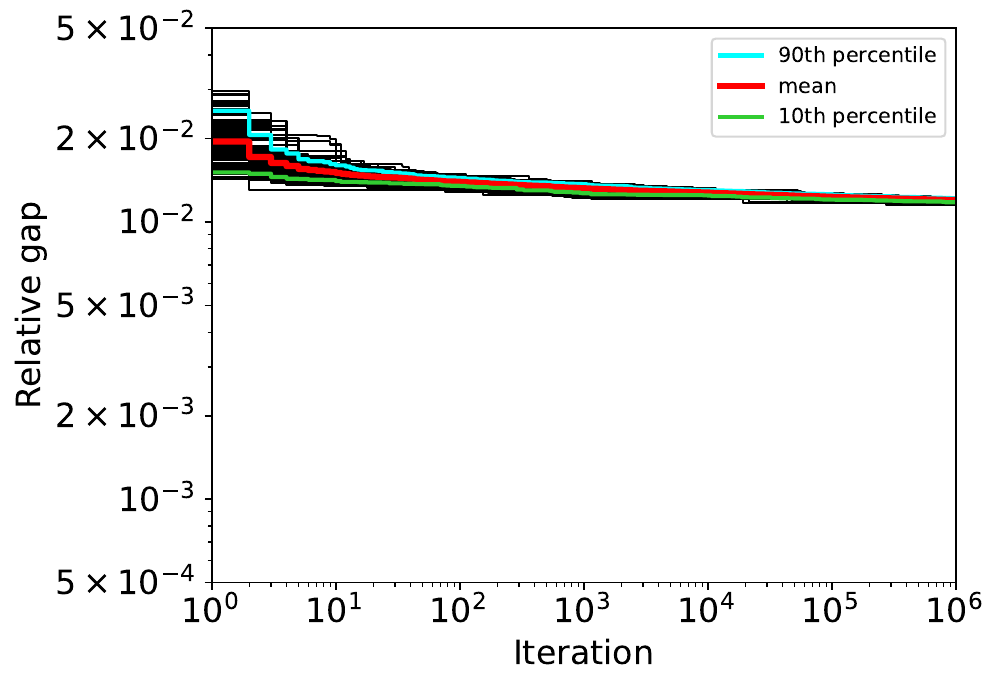}
  \end{minipage}    
  &
  \begin{minipage}{0.42\hsize}
   \centering
   \includegraphics[width=5.5cm, bb=0 0 444 317]{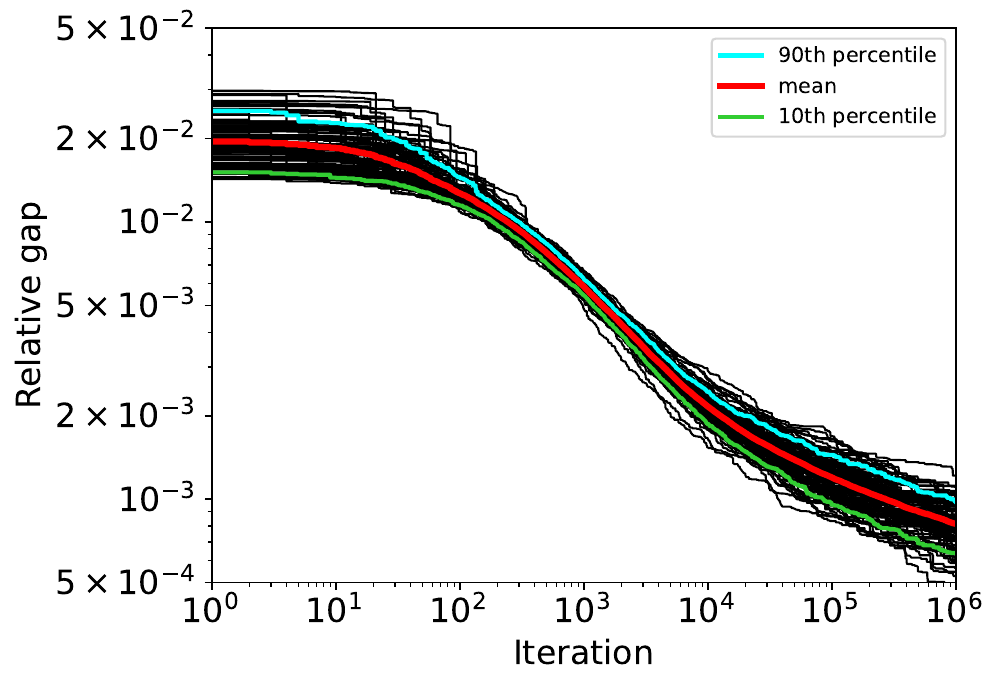}
  \end{minipage}    
  \\
  \hline
  \verb|d18512|
  &
  \begin{minipage}{0.42\hsize}
   \centering
   \includegraphics[width=5.5cm, bb=0 0 444 317]{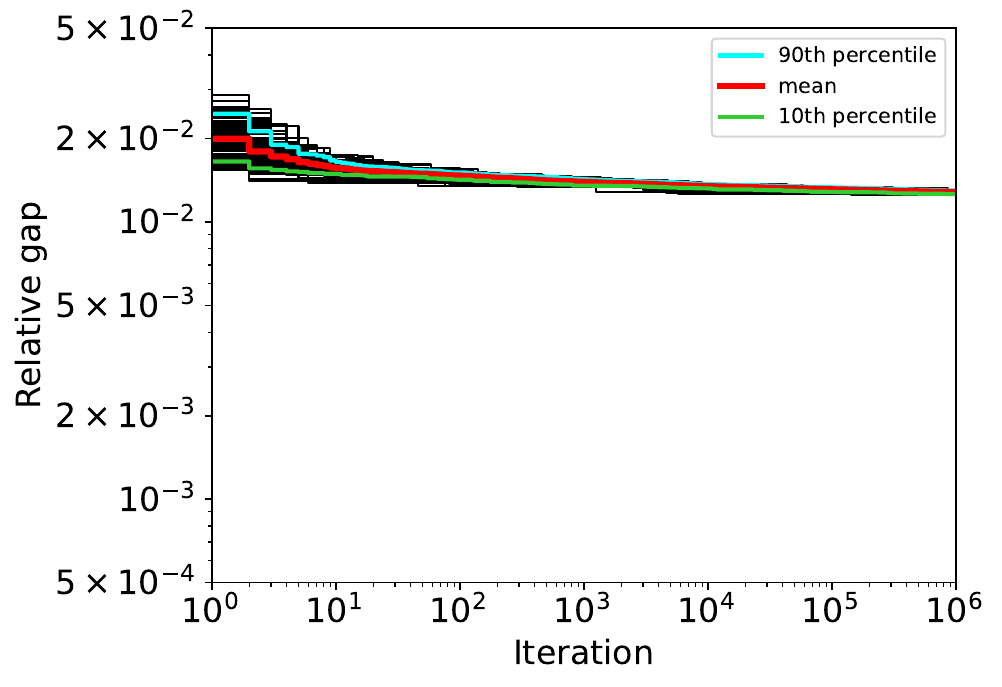}
  \end{minipage}    
  &
  \begin{minipage}{0.42\hsize}
   \centering
   \includegraphics[width=5.5cm, bb=0 0 444 317]{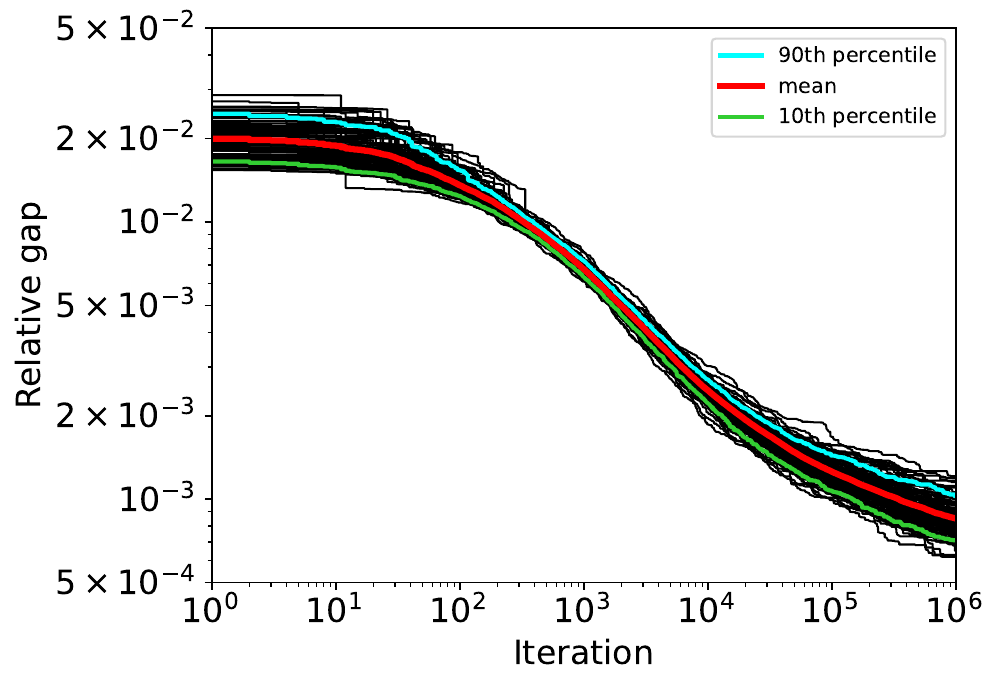}
  \end{minipage}    
  \\
  \hline
  \verb|rl11849|
  &
  \begin{minipage}{0.42\hsize}
   \centering
   \includegraphics[width=5.5cm, bb=0 0 444 317]{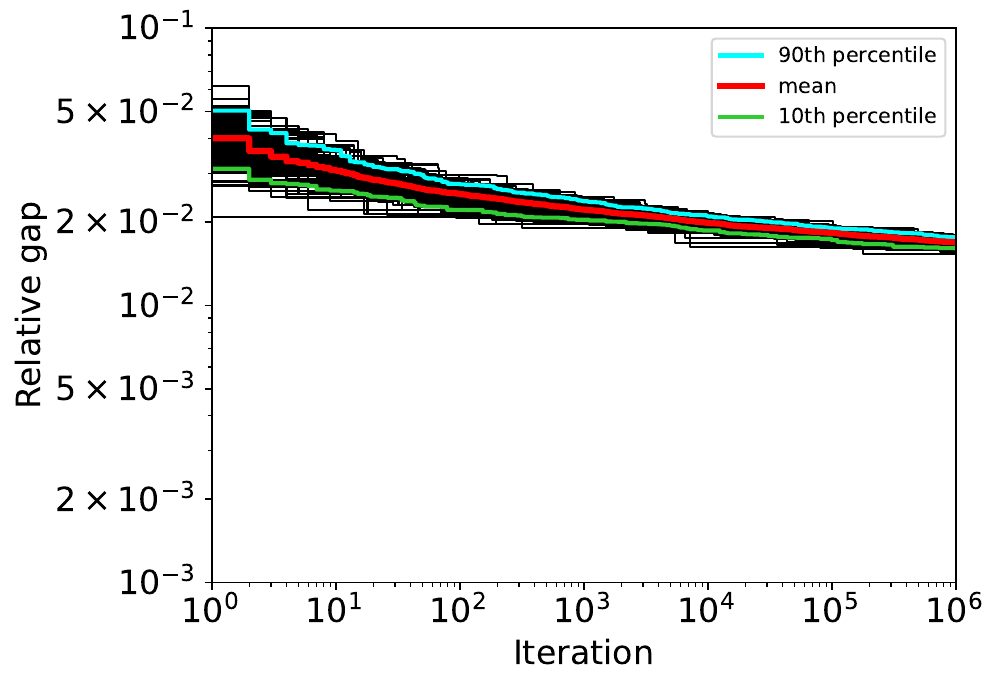}
  \end{minipage}    
  &
  \begin{minipage}{0.42\hsize}
   \centering
   \includegraphics[width=5.5cm, bb=0 0 444 317]{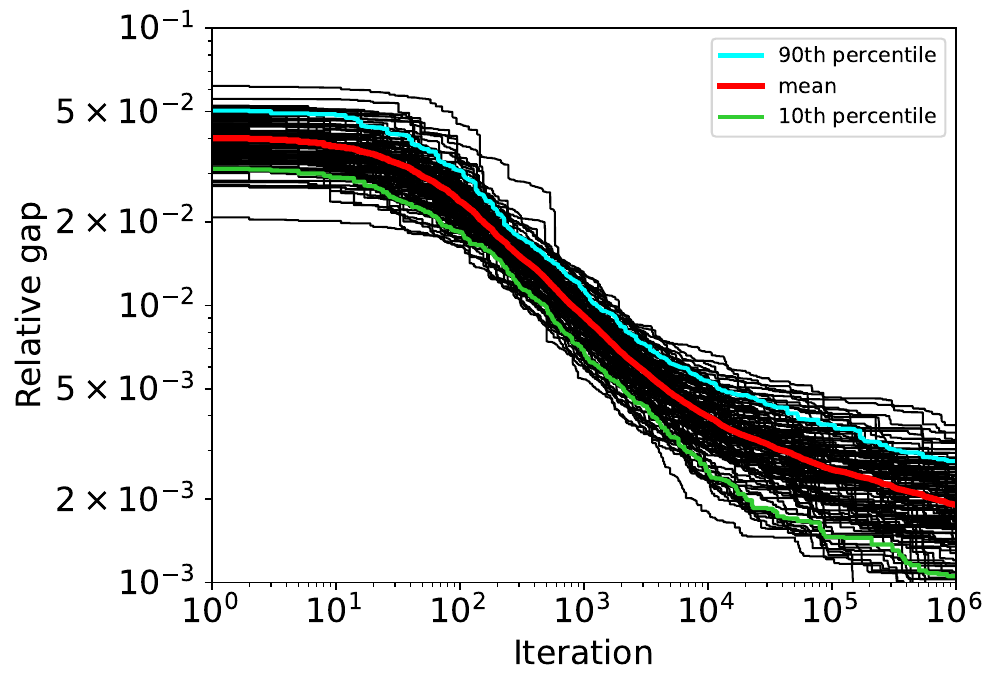}
  \end{minipage}    
  \\
  \hline
  \verb|usa13509|
  &
  \begin{minipage}{0.42\hsize}
   \centering
   \includegraphics[width=5.5cm, bb=0 0 444 317]{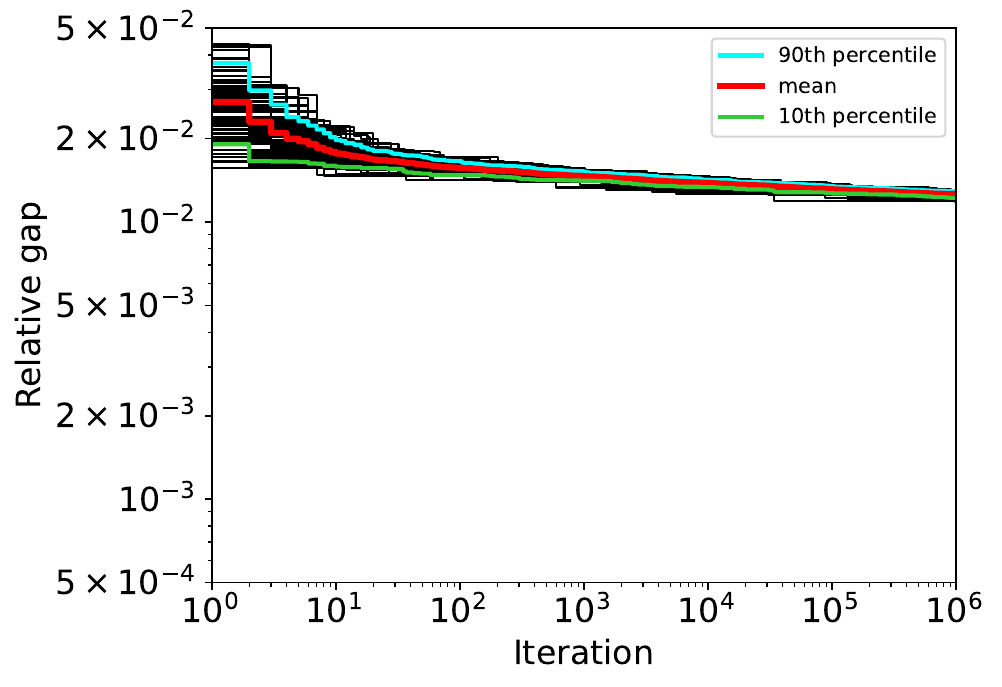}
  \end{minipage}    
  &
  \begin{minipage}{0.42\hsize}
   \centering
   \includegraphics[width=5.5cm, bb=0 0 444 317]{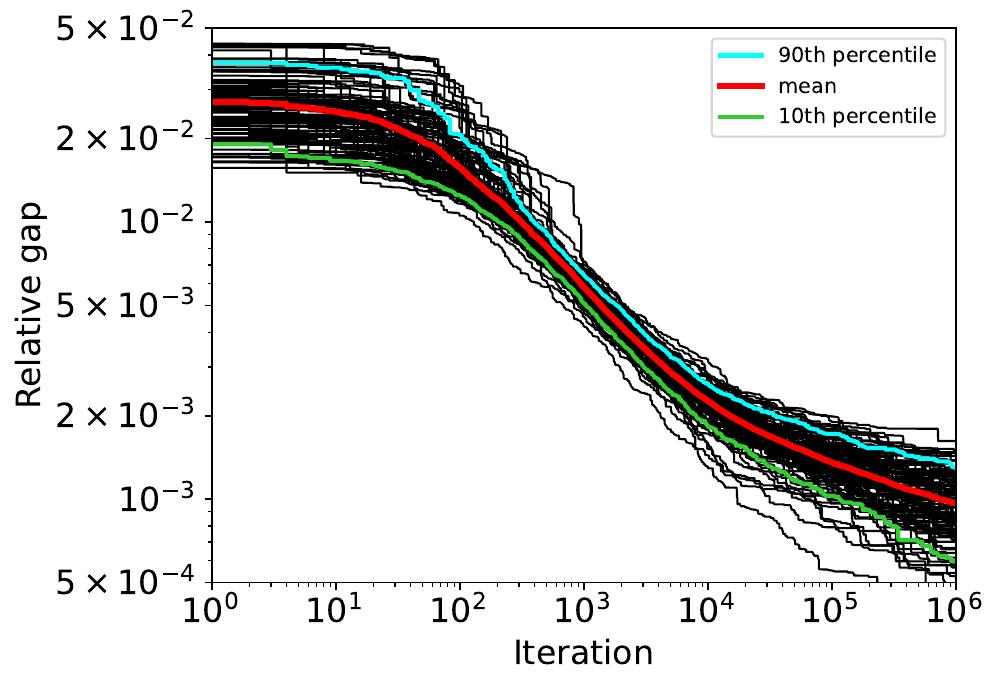}
  \end{minipage}    
  \\
  \hline
 \end{tabular}
 \caption{Evolution of the relative gap of the best EOV over 100 random sets of $10^6$ iterations of the NN+LK RMS and ILS algorithms, respectively, for the five TSPLIB instances listed in Table~\ref{TSPLIB_instances}.}\label{RMS-ILS_NN-LKH}
\end{figure*}

\begin{figure*}[htbp]
 \begin{tabular}{ccc}
  \hline
  & RMS & ILS \\
  \hline
  \verb|brd14051|
  &
  \begin{minipage}{0.42\hsize}
   \centering
   \includegraphics[width=5.5cm, bb=0 0 444 317]{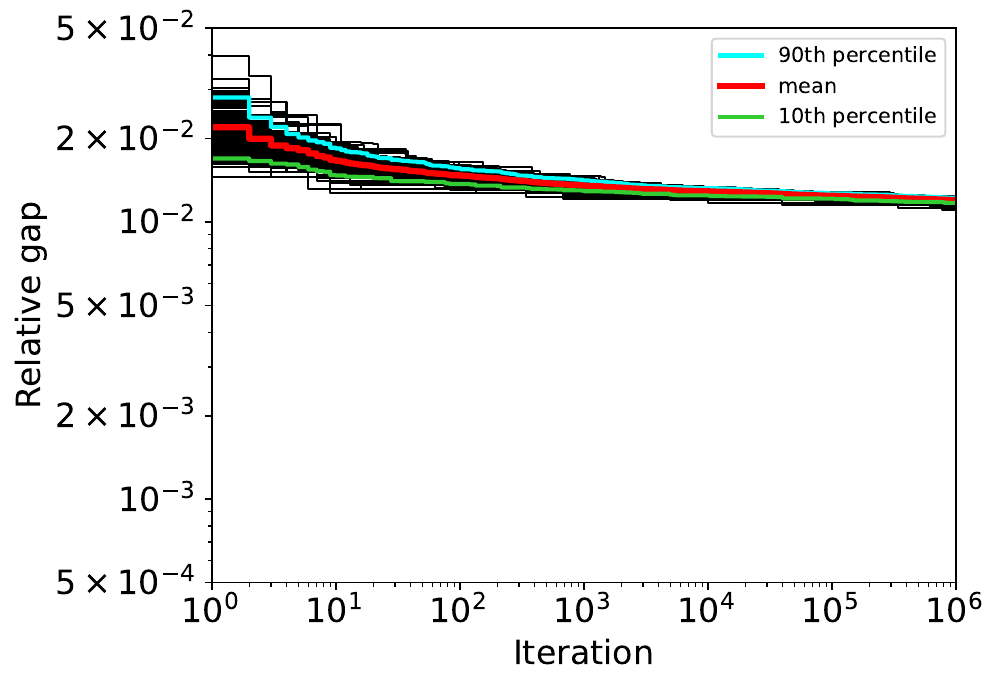}
  \end{minipage}    
  &
  \begin{minipage}{0.42\hsize}
   \centering
   \includegraphics[width=5.5cm, bb=0 0 444 317]{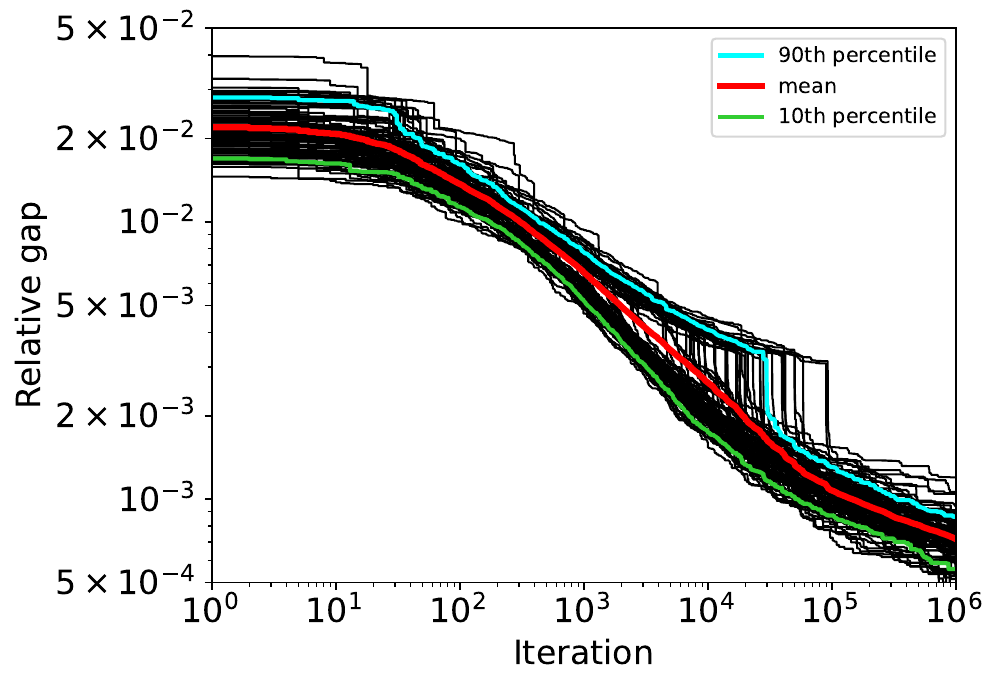}
  \end{minipage}    
  \\
  \hline
  \verb|d15112|
  &
  \begin{minipage}{0.42\hsize}
   \centering
   \includegraphics[width=5.5cm, bb=0 0 444 317]{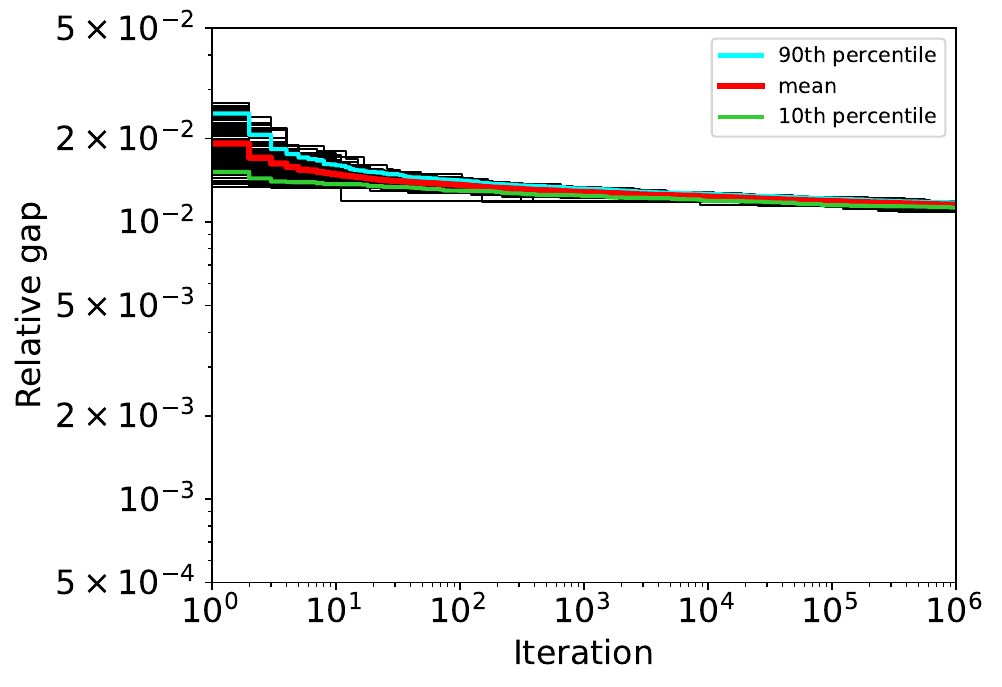}
  \end{minipage}    
  &
  \begin{minipage}{0.42\hsize}
   \centering
   \includegraphics[width=5.5cm, bb=0 0 444 317]{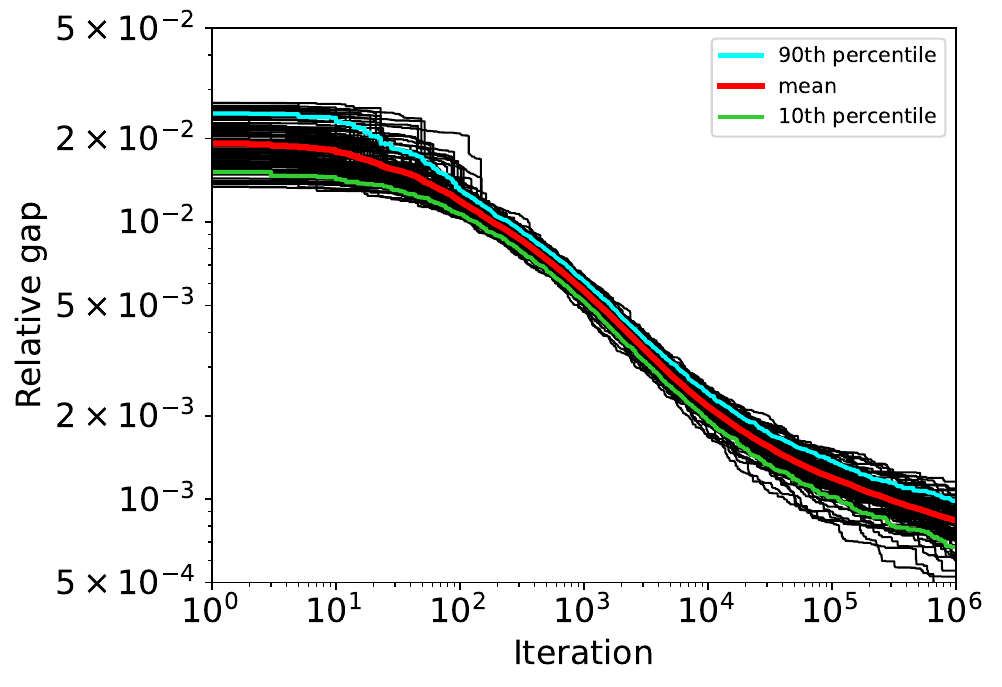}
  \end{minipage}    
  \\
  \hline
  \verb|d18512|
  &
  \begin{minipage}{0.42\hsize}
   \centering
   \includegraphics[width=5.5cm, bb=0 0 444 317]{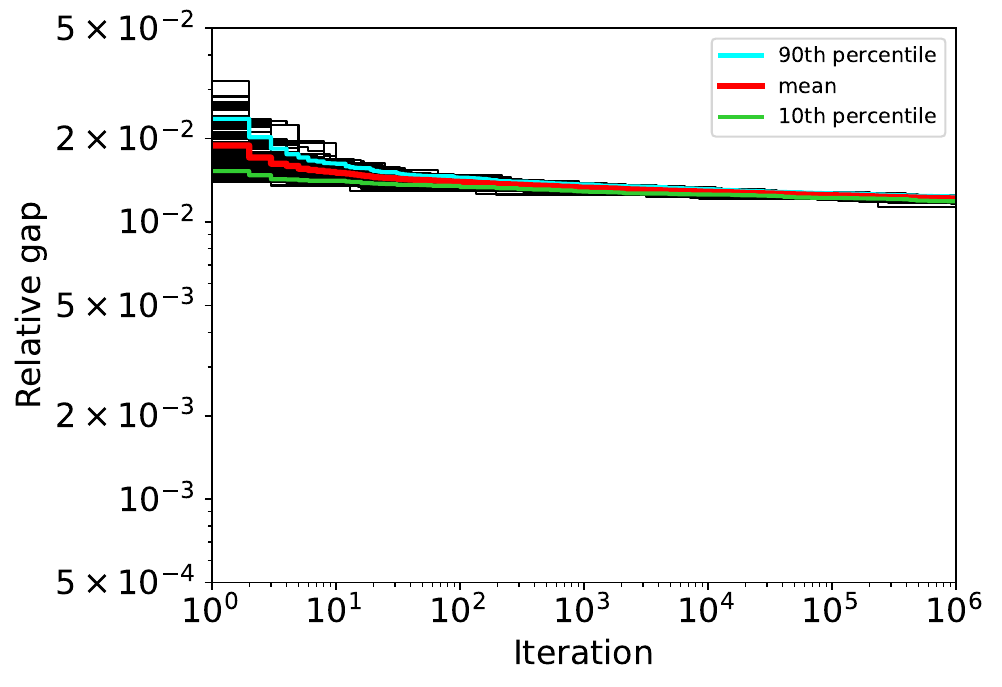}
  \end{minipage}    
  &
  \begin{minipage}{0.42\hsize}
   \centering
   \includegraphics[width=5.5cm, bb=0 0 444 317]{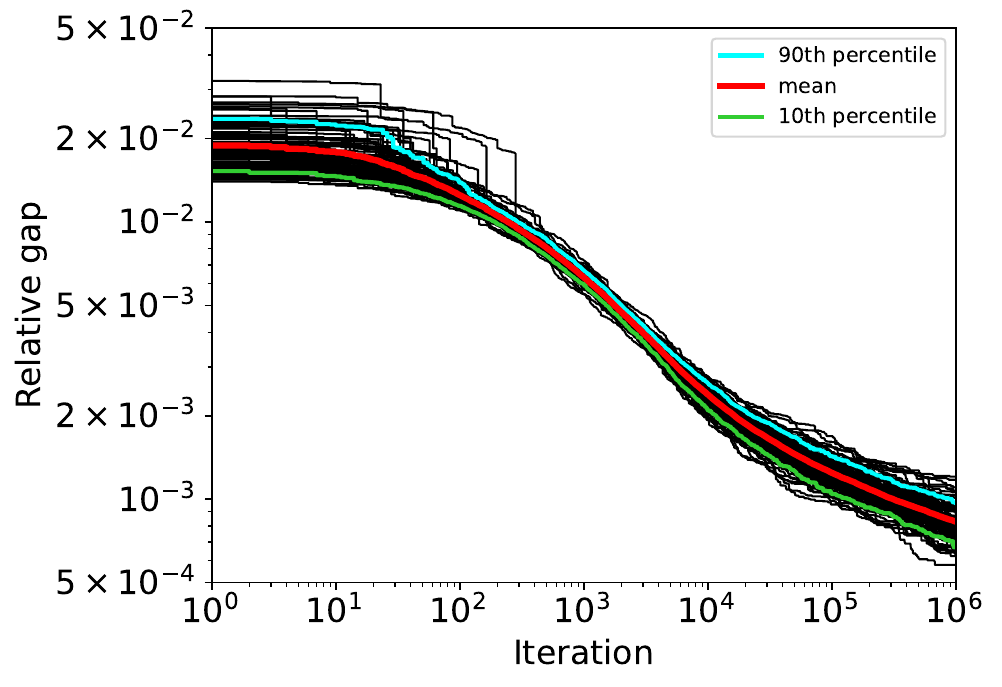}
  \end{minipage}    
  \\
  \hline
  \verb|rl11849|
  &
  \begin{minipage}{0.42\hsize}
   \centering
   \includegraphics[width=5.5cm, bb=0 0 444 317]{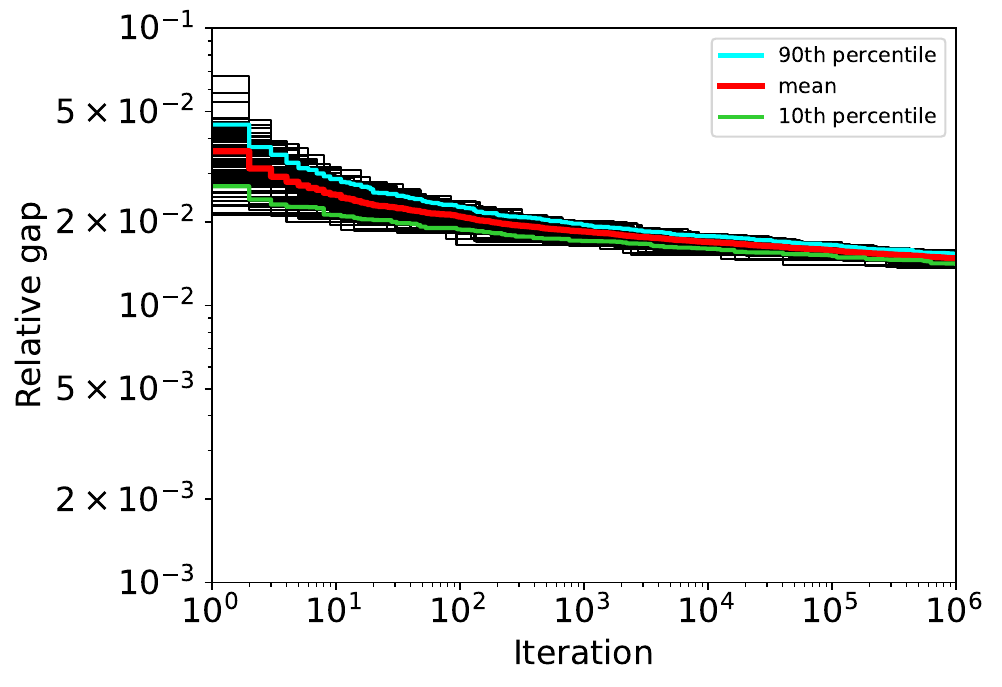}
  \end{minipage}    
  &
  \begin{minipage}{0.42\hsize}
   \centering
   \includegraphics[width=5.5cm, bb=0 0 444 317]{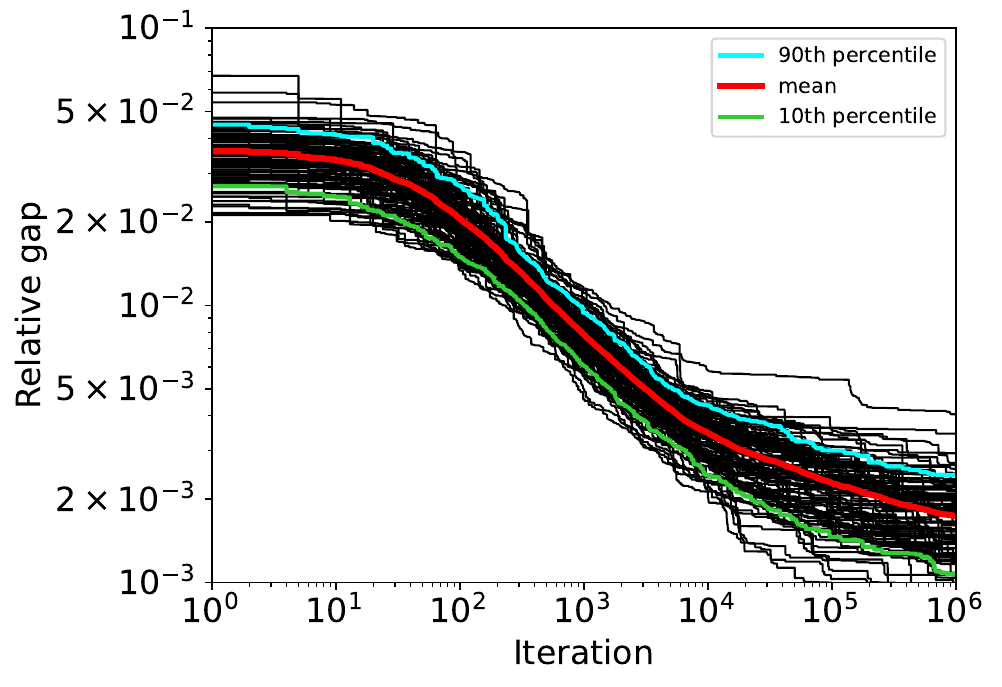}
  \end{minipage}    
  \\
  \hline
  \verb|usa13509|
  &
  \begin{minipage}{0.42\hsize}
   \centering
   \includegraphics[width=5.5cm, bb=0 0 444 317]{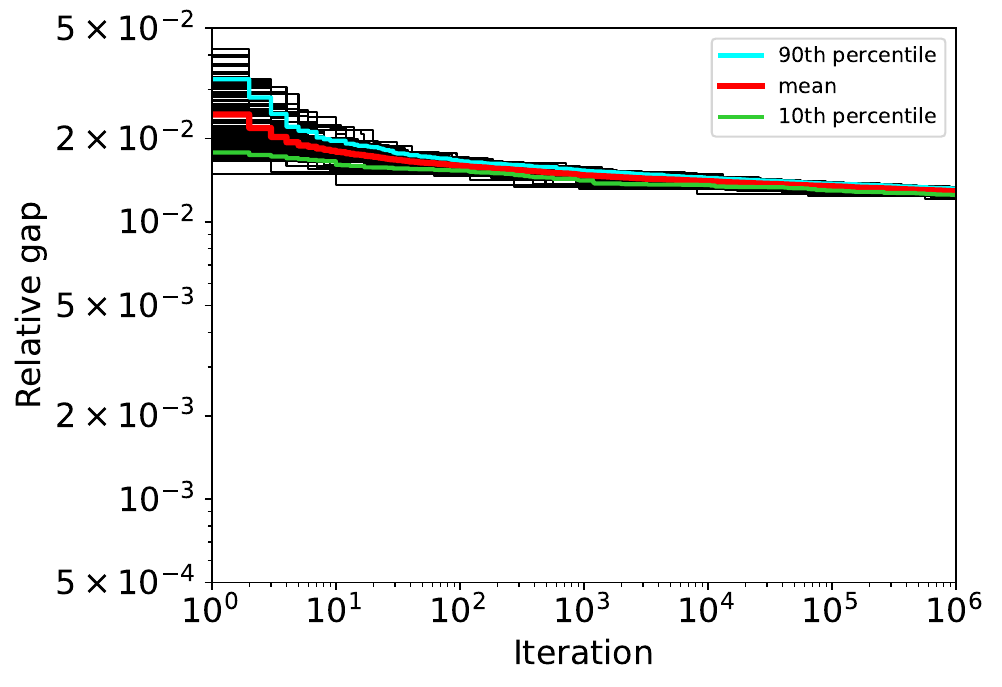}
  \end{minipage}    
  &
  \begin{minipage}{0.42\hsize}
   \centering
   \includegraphics[width=5.5cm, bb=0 0 444 317]{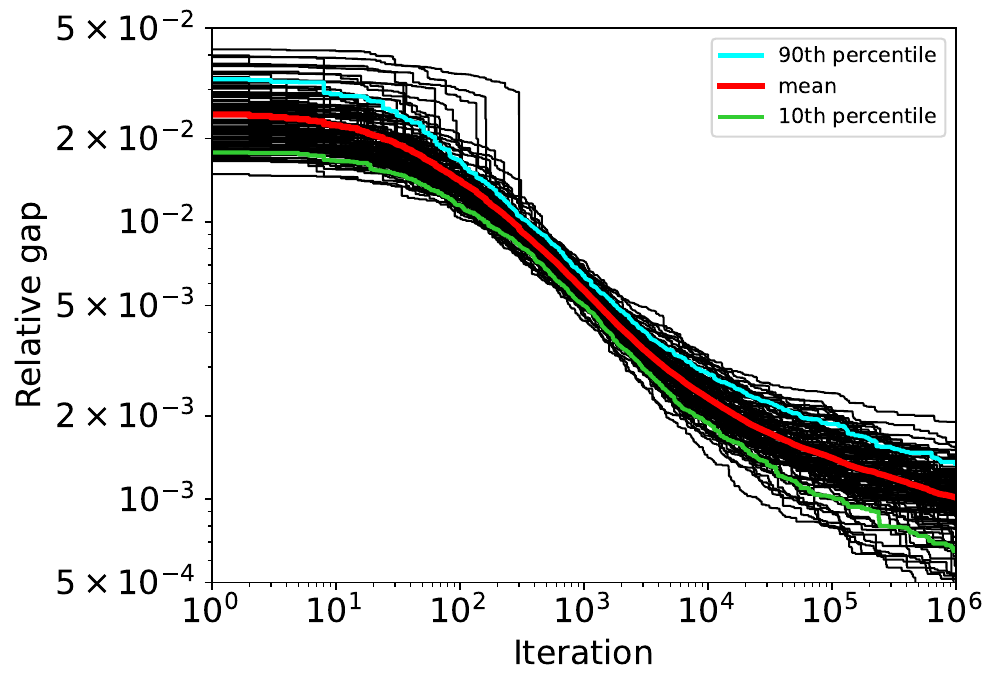}
  \end{minipage}    
  \\
  \hline
 \end{tabular}
 \caption{
Evolution of the relative gap of the best EOV over 100 random sets of $10^6$ iterations of the GR+LK RMS and ILS algorithms, respectively, for the five TSPLIB instances listed in Table~\ref{TSPLIB_instances}.}\label{RMS-ILS_GR-LKH}
\end{figure*}

We perform numerical experiments using a chained Lin-Kernighan (CLK) method \citep[Section~8.3]{Hoos04}, a standard implementation of the ILS method, to examine the intractability of overcoming the curse of scale-freeness. Figures~\ref{RMS-ILS_RA-LKH}--\ref{RMS-ILS_GR-LKH} show the relative gap of the best EOV on a log-log scale when three ILS algorithms and the corresponding three RMS ones are applied to five TSPLIB instances with over 10,000 cities (see Appendix~\ref{subsec_how_to_run_algorithms} for details). Neither method reaches the optimal solution; our RMS algorithms maintain a relative gap above 1\% even after $10^6$ iterations, whereas our ILS algorithms reduce the gap below 0.5\% in all trials and below 0.1\% in some cases. Thus, our ILS algorithms outperform our RMS ones. However, even each of our ILS algorithms exhibits an approximately sigmoidal decay, eventually reducing the gap at best at a power-law rate, which indicates that our ILS algorithms cannot overcome the curse of scale-freeness.

We also investigate why the CLK method reduces the relative gap in an almost sigmoidal manner and fails to overcome the curse of scale-freeness in large-scale optimization. As observed in Figures~\ref{RMS-ILS_RA-LKH}--\ref{RMS-ILS_GR-LKH}, the relative gap decreases very slowly in the early search stage, then declines rapidly (though at best following a power law) in the middle stage, before the decay rate gradually slows. This behavior is due to the depth-first-like exploration of the CLK method, which can be summarized in three steps: (i) by focusing on local regions, the CLK method delays reaching a deep valley where good solutions accumulate because transitions between valleys are not easily made; (ii) once such a deep valley is reached, the best solution is continuously updated; (iii) as unreached good solutions diminish, updates become less frequent. Thus, we infer that the CLK method does not overcome the curse of scale-freeness.

Our theoretical and numerical results indicate that overcoming the curse of scale-freeness is challenging. Our ILS algorithms, based on the CLK method, do not overcome this curse and thus lack effective restart and diversification strategies that exponentially accelerate solution improvement relative to the RMS method. Besides the CLK method, there are several ILS methods for TSP (\citealt{Hoos04}), but it remains unclear whether they overcome the curse. \citet{Dubo15} and \citet{Mu18} reported that state-of-the-art local search algorithms (EAX and LKH) exhibit root-exponential scaling with instance size---with much smaller scaling factors than the exact solver Concorde (\citealt{Cook03}), indicating better scalability---but their results are not directly related to whether these algorithms overcome the curse of scale-freeness.

Then, how can the curse of scale-freeness be overcome? Would it suffice to develop a very powerful LS algorithm with effective restart and diversification strategies? No, it is not that simple. The curse of scale-freeness refers to the phenomenon in which the expected relative gap of the best empirical solution decreases at most in a power-law manner, and our power-law formulas mathematically formalize the conditions under which this ``curse" appears. According to these formulas, even when an algorithm successfully identifies incremental improvements---i.e., even when it efficiently explores the solution space---if the problem is large-scale and the number of potentially generated ``good solutions" (typically local optima) is nearly infinite, the algorithm may still fall victim to the curse of scale-freeness. In other words, the curse comes from the relative difficulty of the problem for the algorithm rather than from a flaw in the algorithm itself. Conversely, if a problem is small or has a favorable structure---thus being relatively easy for the algorithm---the curse of scale-freeness does not appear, which means that such problems can be solved.

\section{Concluding remarks}\label{sec_remarks}

We have introduced the concept of the {\it curse of scale-freeness}, based on the power-law formulas for the expected improvement rate and the expected relative gap of the best empirical objective value (EOV) generated by the RMS method. The curse of scale-freeness can be interpreted as a Zeno's paradox-like phenomenon expressed by the metaphor ``{\it Reaching for the goal makes it slip away}." This concept highlights the inherent intractability of large-scale optimization with multi-start methods.

The curse of scale-freeness suggests that \emph{difficult problems are truly difficult, and only easy problems are easily solved}. Thus, this curse reminds us of the fundamentals of developing metaheuristics for large-scale optimization: adopting a practical strategy to generate the best possible solution within a reasonable computational time. In other words, we should focus on ``when to stop computing" rather than on ``how to close the gap." Moreover, this decision should be guided by a cost-benefit trade-off.

However, this practical strategy cannot be fully realized by simply measuring the gap from the optimal value. As our theoretical and numerical results show, even when there is only a small gap between the optimal value and the best EOV, the computational cost of closing this gap can often exceed expectations. Therefore, the decision to stop an algorithm should be based on the expected improvement rate from additional trials. From this perspective, a promising future task is to develop a method for estimating the expected improvement rate.

\FloatBarrier

%%%%%%%%%%%%%%%%%%%%%%%%%%%%%% Figure 4 %%%%%%%%%%%%%%%%%%%%%%%%%%%%%%
\begin{figure*}[htbp]
\begin{tabular}{cc}
\begin{minipage}[t]{0.48\hsize}
\centering
\includegraphics[scale=0.45,bb=0 0 444 316]{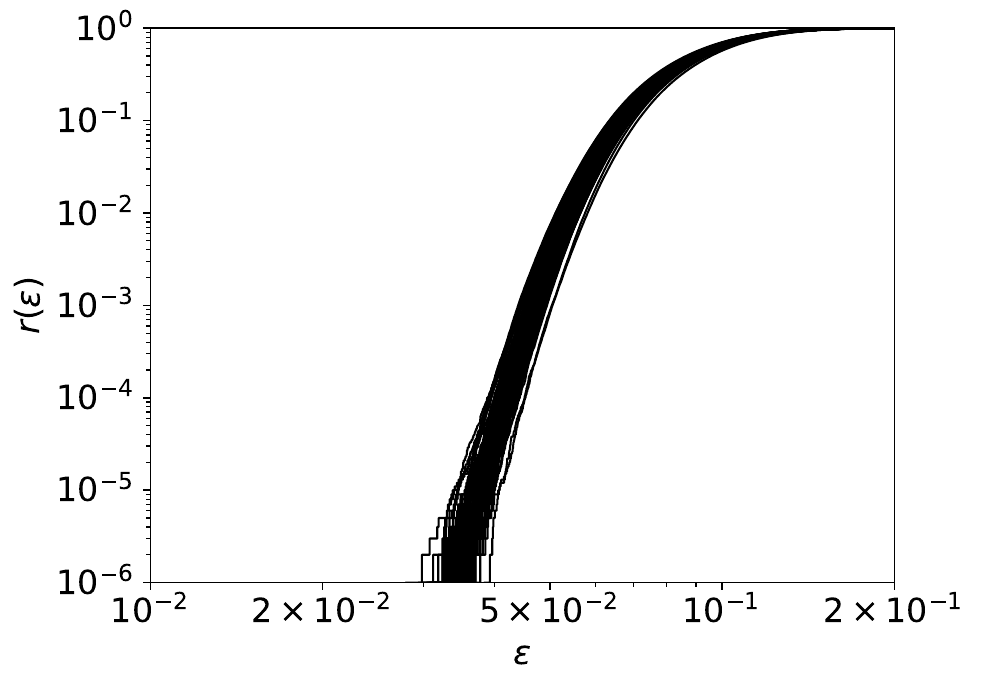}
\subcaption{RA + 3-opt}
\label{fig_good_sol_ratio-RA-3opt}
\end{minipage} &
\begin{minipage}[t]{0.48\hsize}
\centering
\includegraphics[scale=0.45,bb=0 0 444 316]{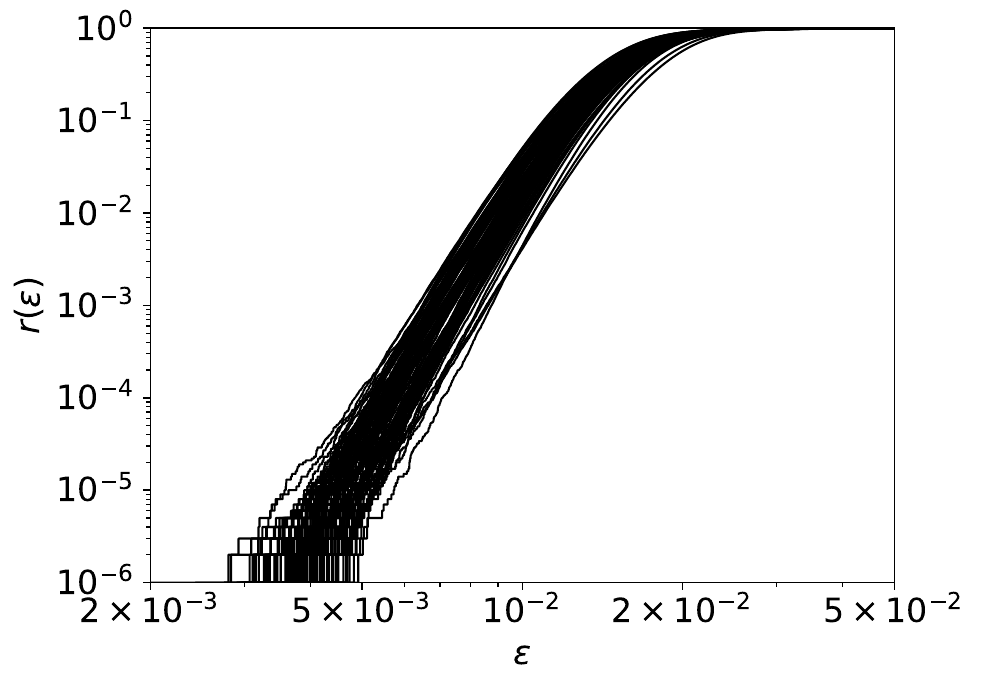}
\subcaption{RA + LK}
\label{fig_good_sol_ratio-RA-LKH}
\end{minipage} \\[10mm]
\begin{minipage}[t]{0.48\hsize}
\centering
\includegraphics[scale=0.45,bb=0 0 444 316]{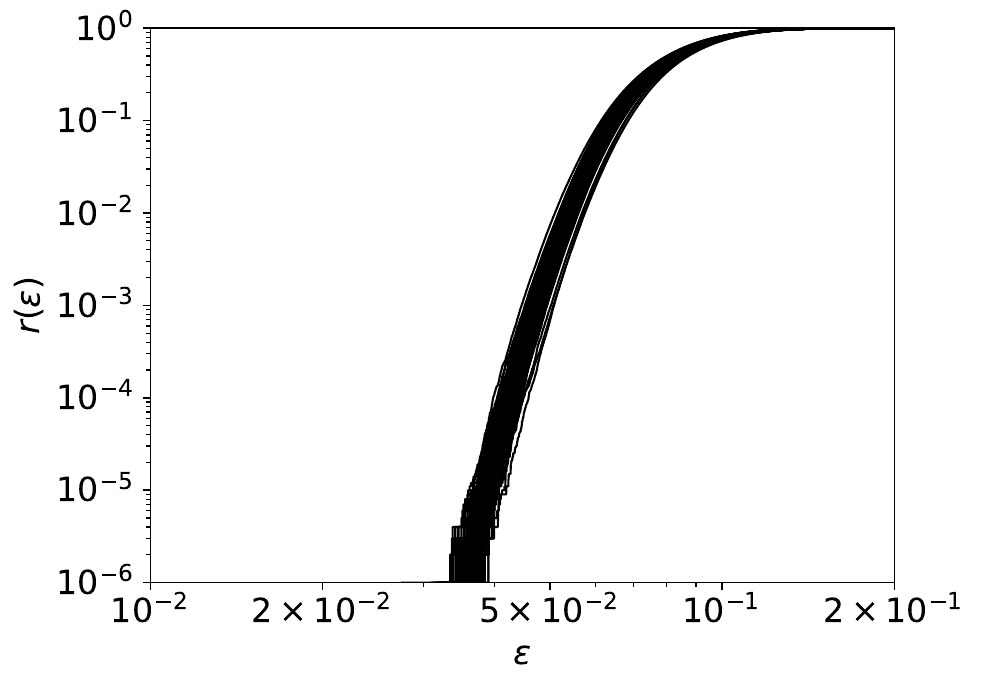}
\subcaption{NN + 3-opt}
\label{fig_good_sol_ratio-NN-3opt}
\end{minipage} &
\begin{minipage}[t]{0.48\hsize}
\centering
\includegraphics[scale=0.45,bb=0 0 444 316]{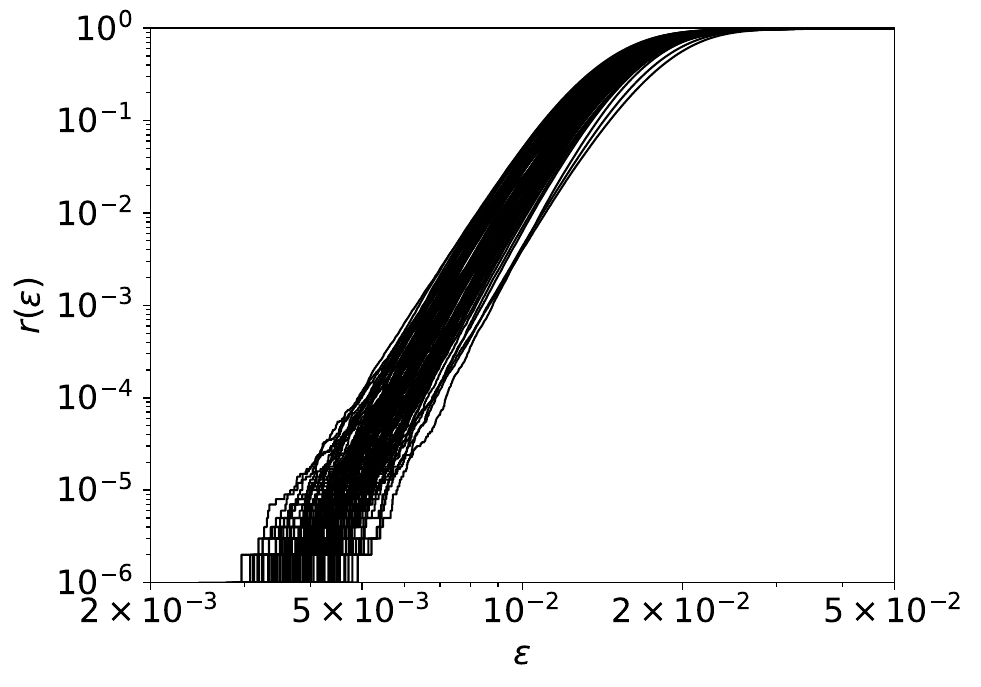}
\subcaption{NN + LK}
\label{fig_good_sol_ratio-NN-LKH}
\end{minipage} \\[10mm]
\begin{minipage}[t]{0.48\hsize}
\centering
\includegraphics[scale=0.45,bb=0 0 444 316]{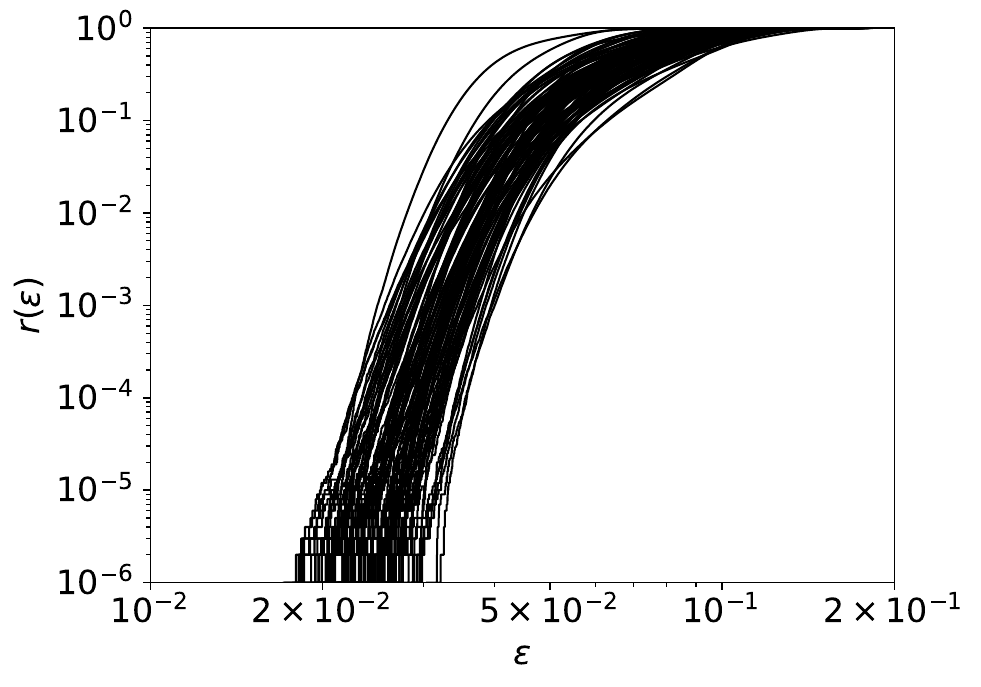}
\subcaption{GR + 3-opt}
\label{fig_good_sol_ratio-GR-3opt}
\end{minipage} &
\begin{minipage}[t]{0.48\hsize}
\centering
\includegraphics[scale=0.45,bb=0 0 444 316]{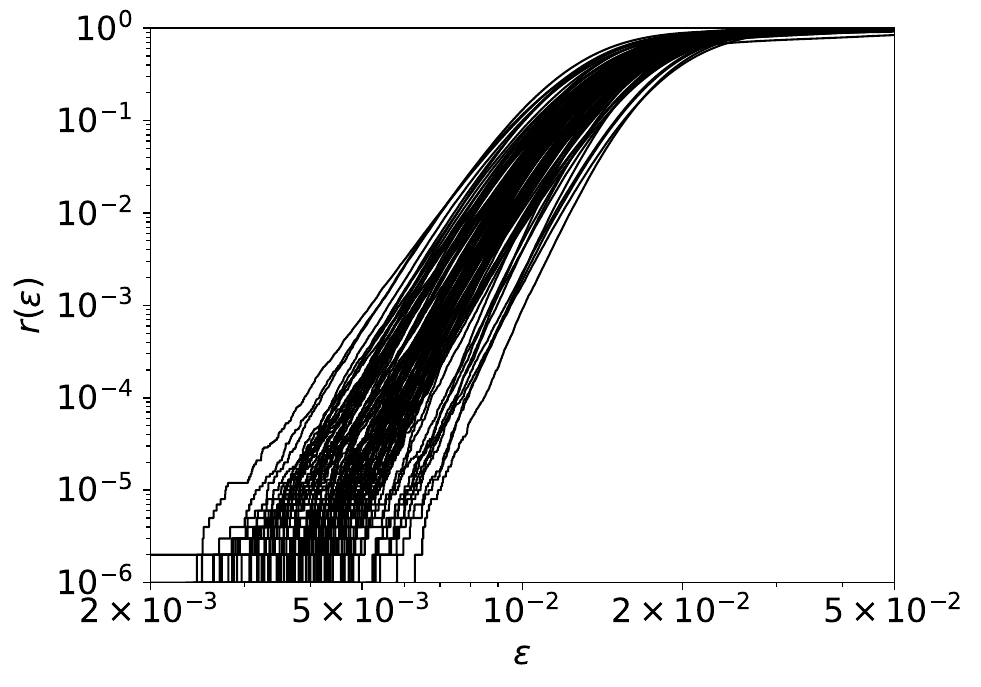}
\subcaption{GR + LK}
\label{fig_good_sol_ratio-GR-LKH}
\end{minipage}
\end{tabular}
\caption{Good-solution-ratio function $r(\varep)$ calculated from the data used to generate Figure~\ref{fig_RMS_1000cities}.}\label{fig_GSR_RMS_1000cities}
\end{figure*}
%%%%%%%%%%%%%%%%%%%%%%%%%%%%%%%%%%%%%%%%%%%%%%%%%%%%%%%%%%%%%%%%%%%%%%

%%%%%%%%%%%%%%%%%%%%%%%%%%%%%% Figure 5 %%%%%%%%%%%%%%%%%%%%%%%%%%%%%%
%
%%%%%%%%%%%%%%%%%%%%%%%%%%%%%%%%%%%%%%%%%%%%%%%%%%%%%%%%%%%%%%%%%%%%%%

%%%%%%%%%%%%%%%%%%%%%%%%%%%%%% Figure 6 %%%%%%%%%%%%%%%%%%%%%%%%%%%%%%
%%%%%%%%%%%%%%%%%%%%%%%%%%%%%%%%%%%%%%%%%%%%%%%%%%%%%%%%%%%%%%%%%%%%%%

%%%%%%%%%%%%%%%%%%%%%%%%%%%%%% Figure 7 %%%%%%%%%%%%%%%%%%%%%%%%%%%%%%

\if0
%%%%%%%%%%%%%%%%%%%%%%%%%%%%%%%%%%%%%%%%%%%%%%%%%%%%%%%%%%%%
\bmhead{Acknowledgments}
Acknowledgments are not compulsory. Where included they should be brief. Grant or contribution numbers may be acknowledged.
%%%%%%%%%%%%%%%%%%%%%%%%%%%%%%%%%%%%%%%%%%%%%%%%%%%%%%%%%%%%
\fi

\appendix

\newpage
\section{Proof of Theorem~\ref{thm_asymp-E[R_n(x)]-01}}\label{proof_thm_asymp-E[R_n(x)]-01}

The proof of Theorem~\ref{thm_asymp-E[R_n(x)]-01} requires the following lemma.
\begin{lem}\label{lem_asymp-E[R_n(x)]}
Suppose that Assumption~\ref{assumpt_GEV} holds, and that $\EE[X] =\int_0^{\infty} x \,d F(x) < \infty$ and $\xi < 1$. For $x \in (0, x^*)$, we then have
\begin{subequations}
\begin{align}
&
\lim_{n\to\infty}
\EE\left[ 
{1 \over a(n)}
\left\{
(Z_n - x)_+ - (U(n) - x)
\right\}
\right]
= \ol{m}_{\xi},
\label{lim_E[(Z_n-x)_+]}
\end{align}
or equivalently,
\begin{align}
\EE[ (Z_n - x)_+ ]
= ( \ol{m}_{\xi} + o(1) ) a(n) + U(n) - x, 
\qquad
\mbox{as $n \to \infty$},
\label{asymp_E[(Z_n-x)_+]}
\end{align}
\end{subequations}
where $\ol{m}_{\xi}$ is given in (\ref{eqn_ol{m}_{xi}}) and $a(\,\cdot\,)$ is a positive function satisfying (\ref{connection-a-to-U}).
\end{lem}

\medskip

\noindent
{\it Proof of Lemma~\ref{lem_asymp-E[R_n(x)]}.} 
We only prove (\ref{lim_E[(Z_n-x)_+]}) since it is equivalent to (\ref{asymp_E[(Z_n-x)_+]}). By splitting the expectation on the left-hand side of (\ref{lim_E[(Z_n-x)_+]}) according to the events $\{Z_n > x\}$ and $\{Z_n \le x\}$, we obtain
\begin{align}
&
\EE\left[ 
{1 \over a(n)} \left\{ (Z_n - x)_+ - (U(n) - x) \right\}
\right]
\notag
\\
&\quad =\EE\left[ \one(Z_n > x) {Z_n - U(n) \over a(n)} \right]
+
\EE\left[ \one(Z_n \le x)
{x - U(n) \over a(n)}
\right]
\notag
\\
&\quad =
\EE\left[ {Z_n - U(n) \over a(n)} \right]
-
\EE\left[ \one(Z_n \le x) {Z_n - U(n) \over a(n)} \right]
+
\EE\left[ \one(Z_n \le x) {x - U(n) \over a(n)}
\right]
\notag
\\
&\quad =
\EE\left[ {Z_n - U(n) \over a(n)} \right]
+
\EE\left[ \one(Z_n \le x) {x - Z_n\over a(n)} \right],\qquad 0 < x < x^*.
\label{eqn_240810-01}
\end{align}
Taking the limit as $n \to \infty$ in the second term on the left-hand side of (\ref{eqn_240810-01}), we obtain the following: For $x \in (0, x^*)$,
\begin{align}
0 \le
\lim_{n\to\infty}
\EE\left[ \one(Z_n \le x) {x - Z_n\over a(n)} \right]
&\le
\lim_{n\to\infty}{x\EE[\one(Z_n \le x)] \over a(n)}
= \lim_{n\to\infty}{x\PP(Z_n \le x) \over a(n)}
\notag
\\
&= \lim_{n\to\infty}{x [ F(x) ]^n \over a(n)} = 0,
\label{eqn_240810-02}
\end{align}
where the last equality holds because the function $a(\,\cdot\,)$ exhibits at most polynomial growth for any $\xi < 1$ (see Propositions~\ref{prop_class_R} and \ref{prop_function_U}). Applying (\ref{eqn_240810-02}) and \cite[Theorem~5.3.1]{Haan06} to (\ref{eqn_240810-01}) yields
\begin{align*}
&
\lim_{n\to\infty}
\EE\left[ 
{1 \over a(n)} \left\{ (Z_n - x)_+ - (U(n) - x) \right\}
\right]
\notag
\\
&\quad =\lim_{n\to\infty}
\EE\left[ {Z_n - U(n) \over a(n)} \right]
= \int_{-\infty}^{\infty}z \,d G_{\xi}(z) = \ol{m}_{\xi},\qquad 0 < x < x^*,
\end{align*}
which shows that (\ref{lim_E[(Z_n-x)_+]}) holds. \qed

\medskip

Using Lemma~\ref{lem_asymp-E[R_n(x)]}, we prove the statements (i)--(iv) of Theorem~\ref{thm_asymp-E[R_n(x)]-01}. To achieve this, we use (\ref{eqn_E[R_n(x)]}), which is derived from (\ref{defn_R_n(x)}), (\ref{eqn_ol{m}_{xi}}), and (\ref{asymp_E[(Z_n-x)_+]}):
\begin{subequations}\label{eqn_E[R_n(x)]}
\begin{empheq}[left = {\EE[R_n(x)] = \empheqlbrace \,}]{alignat = 2}
& 
{1 \over x}
\left[
\left(-\varGamma(-\xi) - {1 \over \xi} + o(1) \right) a(n) + U(n) - x
\right], 
&\quad &   \mbox{$\xi \neq 0,~\xi < 1$},  
\label{eqn_E[R_n(x)]-a} 
\\
& {1 \over x}
\left[
\left(\gamma + o(1) \right) a(n) + U(n) - x
\right],         
&\quad &   \mbox{$\xi = 0$}, 
\label{eqn_E[R_n(x)]-b}
\end{empheq}
\end{subequations}
as $n \to \infty$ for any $x \in (0,x^*)$.

\medskip

\noindent
{\it Proof of Theorem~\ref{thm_asymp-E[R_n(x)]-01}(i).}  Suppose that $0 < \xi < 1$. It then follows from Proposition~\ref{prop_function_U}(ii) that
\begin{align}
\lim_{n\to\infty}
{U(n) - x \over a(n)} = 
\lim_{n\to\infty}
{U(n) \over a(n)} ={1 \over \xi},\qquad  x > 0.
\label{eqn_240224-03}
\end{align} 
Using (\ref{eqn_240224-03}) and (\ref{eqn_E[R_n(x)]-a}), we obtain
\begin{align*}
\lim_{n\to\infty}{1 \over a(n)} \EE[R_n(x)]
&= {1 \over x }
\lim_{n\to\infty}
\left[
-\varGamma(-\xi) - {1 \over \xi} + o(1) + {U(n) - x \over a(n)}
\right]
\nonumber
\\
&= {1 \over x }
\left[
-\varGamma(-\xi) - {1 \over \xi} + {1 \over \xi}
\right]
= -{\varGamma(-\xi) \over x },\qquad x > 0,
\end{align*} 
which shows that (\ref{lim_E[R_n(x)]_xi>0}) holds. \qed

\medskip

\noindent
{\it Proof of Theorem~\ref{thm_asymp-E[R_n(x)]-01}(ii).}
Suppose that $\xi = 0$ and $x^* = \infty$. It then follows from (\ref{eqn_lim_U(t)=x^*}) and (\ref{eqn_240224-01}) (see Proposition~\ref{prop_function_U}(iv)) that
\begin{align}
\lim_{n\to\infty}U(n) = x^* = \infty,
\qquad
\lim_{n\to\infty} {a(n) \over U(n)} = 0.
\label{eqn_240224-04}
\end{align} 
Using (\ref{eqn_240224-04}) and (\ref{eqn_E[R_n(x)]-b}), we obtain
\begin{align*}
\lim_{n\to\infty}{1 \over U(n)} \EE[R_n(x)]
&=
{1 \over x}
\lim_{n\to\infty}
\left[
(\gamma + o(1)) {a(n) \over U(n)} + 1 - {x \over U(n)}
\right]
= {1 \over x},
\qquad x > 0,
\end{align*} 
which shows that (\ref{lim_E[R_n(x)]_xi=0_infinite}) holds. \qed

\medskip

\noindent
{\it Proof of Theorem~\ref{thm_asymp-E[R_n(x)]-01}(iii).}
Suppose that $\xi = 0$ and $x^* < \infty$.  It then follows from (\ref{eqn_lim_U(t)=x^*}) and (\ref{eqn_240224-02}) (see Proposition~\ref{prop_function_U}(iv)) that
\begin{align}
\lim_{n\to\infty} {a(n) \over x^* - U(n)} = 0.
\label{eqn_240224-05}
\end{align} 
Using (\ref{eqn_240224-05}) and (\ref{eqn_E[R_n(x)]-b}), we obtain
\begin{align*}
&
\lim_{n\to\infty}
{1  \over x^* - U(n)} \left[ {x^* - x \over x} - \EE[R_n(x)] \right]
\nonumber
\\
&\quad= 
{1 \over x}
\lim_{n\to\infty}
{1  \over x^* - U(n)}
\left[ x^* - x - (\gamma + o(1)) a(n) - U(n) + x \right]
\nonumber
\\
&\quad= 
{1 \over x}
\lim_{n\to\infty}
{1  \over x^* - U(n)}
\left[ x^* - U(n) - \gamma a(n)\right]
\nonumber
\\
&\quad= 
{1 \over x}
\lim_{n\to\infty}
\left[ 1 - { \gamma a(n)\over x^* - U(n)} \right]
=
{1 \over x},\qquad 0 < x < x^*,
\end{align*} 
which shows that (\ref{lim_E[R_n(x)]_xi=0_finite}) holds.
Note here that (\ref{eqn_affine_relation}) yields
\begin{align}
\EE[\varDelta_n(x)] = {x \over x^* }
\left[ {x^* - x \over x} - \EE[R_n(x)] \right],\qquad 0 < x < x^*.
\label{eqn_240224-06}
\end{align} 
Combining (\ref{eqn_240224-06}) and (\ref{lim_E[R_n(x)]_xi=0_finite}) leads to 
(\ref{lim_E[Delta_n(x)]_xi=0_finite}). \qed

\medskip

\noindent
{\it Proof of Theorem~\ref{thm_asymp-E[R_n(x)]-01}(iv).}
Suppose that $\xi < 0$. It then follows from Proposition~\ref{prop_function_U}(iii) that
\begin{align}
\lim_{n\to\infty} U(n) = x^* < \infty,\qquad
\lim_{n\to\infty} {x^* - U(n) \over a(n)} = {-1 \over \xi}.
\label{eqn_240901-01}
\end{align}
Using (\ref{eqn_240901-01}) and (\ref{eqn_E[R_n(x)]-a}), we obtain
\begin{align*}
&
\lim_{n\to\infty}
{1 \over a(n)} \left[ {x^* - x \over x} - \EE[R_n(x)] \right]
\nonumber
\\
&\quad= 
{1 \over x}
\lim_{n\to\infty}
{1 \over a(n)}
\left[ x^ * - x - 
\left(-\varGamma(-\xi) - {1 \over \xi} + o(1) \right)a(n) - U(n) + x  
\right]
\nonumber
\\
&\quad= 
{1 \over x}
\lim_{n\to\infty}
\left[ {x^* - U(n) \over a(n)} + \varGamma(-\xi) + {1 \over \xi} \right]
\nonumber
\\
&\quad= 
{1 \over x}
\left[ -{1 \over \xi} + \varGamma(-\xi) + {1 \over \xi} \right]
={ \varGamma(-\xi) \over x},\qquad 0 < x < x^*,
\end{align*} 
which shows that (\ref{lim_E[R_n(x)]_xi<0}) holds. Furthermore, 
combining (\ref{lim_E[R_n(x)]_xi<0}) with (\ref{eqn_240224-06}) leads to (\ref{lim_E[Delta_n(x)]_xi<0}).  The proof of Theorem~\ref{thm_asymp-E[R_n(x)]-01} is completed. \qed

\section{Proofs of propositions}

\subsection{Proof of Proposition~\ref{prop_scale-free-01}}\label{proof:prop_scale-free-01}

First, we prove (\ref{scale-free-ERG-01}). It follows from (\ref{ean_ERG_in_R_{xi}}) and $L \in \calR_0$ that, for all $c \in \bbN$,
\begin{align*}
\lim_{n\to\infty} {\EE[\varDelta_{cn}(x)] \over \EE[\varDelta_n(x)]}
=  \lim_{n\to\infty} {L(cn) \over L(n)} { (cn)^{\xi} \over n^{\xi} }
= c^{\xi} \lim_{n\to\infty} {L(cn) \over L(n)}
=  c^{\xi},\qquad 0 < x < x^*,
\end{align*}
which shows that (\ref{scale-free-ERG-01}) holds. 

Next, we prove (\ref{scale-free-ERG-02}). It follows from (\ref{ean_ERG_in_R_{xi}}) that
\begin{align}
\log \EE[\varDelta_n(x)] = \log L(n) + \xi \log n. 
\label{lim_log_E[Delta_n(x)]_xi<0}
\end{align}
Since $L \in \calR_0$, $\log L(t)$ can be expressed in the following form, known as the Karamata representation (see \citealt[Theorem~B.1.6]{Haan06}):
\begin{align}
\log L(t) = \log \eta (t) + \int_{t_0}^t {\theta(s) \over s} ds,
\qquad t > t_0,
\label{eqn_240228-01}
\end{align}
for some $t_0 > 0$, where $\eta: \bbR_+ \to \bbR$ and $\theta: \bbR_+ \to \bbR$ are measurable functions such that $\lim_{t \to \infty}\eta(t) =: \eta_0 \in (0,\infty)$ and $\lim_{t \to \infty}\theta(t) =0$ (the latter limit equals the index of regular variation of $L \in \calR_0$). Applying these limits to (\ref{eqn_240228-01}) yields
\begin{align*}
\lim_{t\to\infty}{\log L(t) \over \log t}
&= \lim_{t\to\infty} {\log \eta (t) \over \log t} 
+ \lim_{t\to\infty} {1 \over \log t} \int_{t_0}^t {\theta(s) \over s} ds
\nonumber
\\
&= \lim_{t\to\infty} {1 \over t^{-1}} {\theta(t) \over t}
= \lim_{t\to\infty} \theta(t) = 0,
\end{align*}
where the second equality holds by L'H\^{o}pital's rule. Combining the above equation with (\ref{lim_log_E[Delta_n(x)]_xi<0}) leads to (\ref{scale-free-ERG-02}).

\subsection{Proof of Proposition~\ref{prop_power-law-decay}}\label{proof:prop_power-law-decay}

Using (\ref{defn_r(varep)}), we rewrite (\ref{scale-Free-r(varep)}) as follows:
\begin{align*}
c^{\psi}
=
\lim_{\varep \downarrow 0}
{ 1 - F(x^* - c\varep x^*) 
\over 
1 - F(x^* - \varep x^*) 
}
=\lim_{t \downarrow 0}{1 - F(x^* - ct) 
\over 
1 - F(x^* - t)
},
\qquad c > 0,
\end{align*} 
which implies that $F \in \mathsf{MDA}(G_{\xi})$ with $\xi = -1/\psi < 0$ (\citealt[Theorem~1.2.1]{Haan06}). Thus, Assumption~\ref{assumpt_GEV} holds for $G = G_{\xi}$ with $\xi = -1/\psi < 0$ (see Proposition~\ref{prop_GP}). Therefore, the conditions of Proposition~\ref{prop_power-law-decay} imply those of Proposition~\ref{prop_scale-free-01}, and hence, (\ref{scale-free-ERG-01}) and (\ref{scale-free-ERG-02}) hold with $\xi = -1/\psi < 0$.

\subsection{Proof of Proposition~\ref{prop_exponential-type-decay}}\label{proof:prop_exponential-type-decay}

Substituting (\ref{cond_F_moderate_exp}) into (\ref{eqn:F_expressed_by_r}) yields
\begin{align}
1 - F(x) 
= c_0 \exp \left\{-g\left( {x^* \over x^* - x} \right) \right\},
\qquad 0 \le x < x^*.
\label{eqn_240831-01}
\end{align}
Assume the existence of a function $f: \bbR_+ \to (0,\infty)$ such that 
\begin{align}
1 - F(x) 
= c_0 \exp \left\{-\int_0^x  {d s \over f(s)} \right\},
\qquad 0 \le x < x^*.
\label{eqn_240831-02}
\end{align}
It then follows from (\ref{eqn_240831-01}) and (\ref{eqn_240831-02}) that
\begin{align*}
\int_0^x  {d s \over f(s)} = g\left( {x^* \over x^* - x} \right), 
\qquad 0 \le x < x^*.
\end{align*}
Differentiating both sides of the above equation yields 
\begin{align*}
{1 \over f(x)} = {x^* \over (x^* - x)^2} g'\left( {x^* \over x^* - x} \right), 
\qquad 0 < x < x^*.
\end{align*}
Thus,
\begin{align}
f(x) 
&= {(x^* - x)^2 \over x^*} {1 \over g'\left( \dm{x^* \over x^* - x} \right)}, 
\qquad 0 < x < x^*,
\label{eqn_240831-03}
\\
f'(x) 
&= -{2(x^* - x) \over x^*} {1 \over g'\left( \dm{x^* \over x^* - x} \right)}
- {g''\left( \dm{x^* \over x^* - x} \right) \over \left[ g'\left( \dm{x^* \over x^* - x} \right) \right]^2}, 
\qquad 0 < x < x^*.
\label{eqn_240831-04}
\end{align}
Letting $t = x^* / (x^* - x)$, it follows from (\ref{cond_g(t)-01}), (\ref{cond_g(t)-02}), (\ref{eqn_240831-03}), and (\ref{eqn_240831-04}) that 
\begin{align*}
\lim_{x \uparrow x^*} f(x) 
&= \lim_{t \to \infty}  {x^* \over t^2 g'(t)} = 0,
\\
\lim_{x \uparrow x^*} f'(x) 
&= - \lim_{t \to \infty} {2 \over tg'(t)}
- \lim_{t \to \infty} {g''(t) \over \left[ g'(t) \right]^2}
= 0.
\end{align*}
The combination of these equations and (\ref{eqn_240831-02}) implies that the distribution $F$ is a von Mises function in $\mathsf{MDA}(G_0)$ and thus $F \in \mathsf{MDA}(G_0)$ (see \citealt[Theorem~1.2.6]{Haan06}). Therefore, the conditions of Theorem~\ref{thm_asymp-E[R_n(x)]-01}(iii) are satisfied, and (\ref{lim_E[Delta_n(x)]_xi=0_finite}) yields
\begin{align*}
\EE[\varDelta_n(x)]
&= \left( {1 \over x^* } + o(1) \right) (x^* - U(n)),
\qquad \mbox{as $n \to \infty$},
\end{align*}
for $x \in (0,x^*)$, where $U(t) = x^* - L(t)$ for some $L \in \calR_0$ such that $\lim_{t\to\infty}L(t)=0$. As a result, the statement of this proposition holds.

\section{Setup of numerical experiments}\label{sec_num_setup}

This section describes the setup of our numerical experiments aimed at verifying whether the power-law phenomenon---referred to as the {\it curse of scale-freeness}---occurs when improving the best EOV. First, we introduce the Traveling Salesman Problem (TSP) instances used to evaluate the practical implications of our theoretical results. Next, we present the RMS and ILS algorithms implemented in the Concorde TSP Solver (\citealt{Cook03}), a software package that includes an efficient branch-and-cut algorithm along with several heuristic algorithms for TSP. Finally, we outline the procedure for running our RMS and ILS algorithms. Note that these algorithms are executed to obtain empirical solutions that support the phenomena suggested by our theoretical analysis.

\subsection{TSP instances}\label{subsec_instances}

TSP is one of the most well-known combinatorial optimization problems \citep{Law85,Appl07}. Its difficulty can be easily adjusted, and extensive collections of benchmark instances with known optimal solutions are readily available. Therefore, we chose TSP to validate our mathematical results from a practical perspective. Note that while TSP is originally a minimization problem, EVT---our mathematical tool---is designed for maximum value data. To fit TSP into our EVT framework, we multiplied the objective function by $-1$.

To generate Figures~\ref{fig_power-law}--\ref{fig_GSR_RMS_1000cities}, we created 100 random TSP instances as follows. Each instance consists of 1000 cities with two-dimensional coordinates $(x^{(1)}, x^{(2)})$ that are uniformly distributed random integers chosen from $[0, 10^5-1]$. The travel cost $c_{i,j}$ between cities $i$ and $j$, with coordinates $(x_i^{(1)}, x_i^{(2)})$ and $(x_j^{(1)}, x_j^{(2)})$, respectively, is defined as
\begin{align*}
c_{i,j}
= \left\lfloor 
\sqrt{
\left( x_i^{(1)}- x_j^{(1)} \right)^2 + 
\left( x_i^{(2)}- x_j^{(2)} \right)^2
} + 0.5 
\right\rfloor.
\end{align*}
We computed the optimal values for all 100 instances using the branch-and-cut algorithm of the Concorde TSP Solver \citep{Appl07}.

To generate Figures~\ref{RMS-ILS_RA-LKH}--\ref{RMS-ILS_GR-LKH}, we used the five TSP instances listed in Table~\ref{TSPLIB_instances} from TSPLIB \citep{Rei91}, a well-known and widely used benchmark set of TSP instances. These instances contain more than 10,000 cities, and their optimal values are already known.

\begin{table}[h] 
\caption{TSPLIB instances} 
\centering 
\begin{tabular}{ll} 
\hline 
Instance Name & Number of Cities 
\\
\hline 
brd14051 & 14,051 \\ 
d15112 & 15,112 \\
d18512 & 18,512 \\
rl11849 & 11,849 \\ 
usa13509 & 13,509 
\\ 
\hline 
\end{tabular} 
\label{TSPLIB_instances} 
\end{table}

\subsection{Algorithms}\label{subsec_algorithms}

This subsection describes the algorithms used in our experiments on the TSP. We implemented six RMS algorithms by combining three initial solution generators and two local search methods (3-opt and LK) for the RMS-only experiments. For the RMS vs. ILS comparison, we used the three algorithms (with LK as the local search method) for each method. Table~\ref{tb_RMS_algorithms} summarizes all implemented algorithms for both the RMS-only and RMS vs.\ ILS experiments. All algorithms were obtained by modifying the heuristic procedures provided by the Concorde TSP Solver. The three initial solution generators and two local search methods are described in detail below (see \citealt{John03} for more information).

\begin{table*}[tb]
\caption{The options for initial solution generation and local search in RMS and ILS}
\begin{center}
\begin{tabular}{llll}
\hline
Method & Algorithm name & Initial solution generation & Local search \\ \hline
RMS & RA + 3-opt & RAndom & 3-opt \\
RMS & NN + 3-opt & randomized Nearest Neighbor & 3-opt \\
RMS & GR + 3-opt & randomized GReedy & 3-opt \\
RMS & RA + LK    & Random & Lin-Kernighan \\
RMS & NN + LK    & randomized Nearest Neighbor & Lin-Kernighan \\
RMS & GR + LK    & randomized GReedy & Lin-Kernighan 
\\ \hline\hline
ILS & RA + LK    & RAndom & Lin-Kernighan \\
ILS & NN + LK    & randomized Nearest Neighbor & Lin-Kernighan \\
ILS & GR + LK    & randomized GReedy & Lin-Kernighan \\ \hline
\end{tabular}
\end{center}
\label{tb_RMS_algorithms}
\end{table*}

The three initial solution generators are (i) the RAndom (RA) algorithm, (ii) the randomized Nearest Neighbor (NN) algorithm, and (iii) the randomized GReedy (GR) algorithm. The RA algorithm is originally provided by the Concorde TSP Solver, while the NN and GR algorithms have been modified to incorporate a random mode. The characteristics of these three algorithms are as follows:
\begin{enumerate}
\item The RA algorithm generates an initial solution by visiting each city exactly once in a random order.
\item The NN algorithm starts from a randomly selected city and moves to a randomly selected city among the three nearest unvisited ones.
\item The GR algorithm starts with the shortest edge and adds a randomly selected edge among the three shortest remaining ones to the current path, avoiding the creation of subtours and vertices of degree three.
\end{enumerate}

\subsection{How to run RMS and ILS algorithms}\label{subsec_how_to_run_algorithms}

This subsection describes how the RMS and ILS algorithms were run to obtain the data shown in Figures~\ref{fig_RMS_1000cities}--\ref{RMS-ILS_GR-LKH}, with a focus on the use of random numbers generated by the Mersenne Twister.

We replaced the default random number generator (\citealt[Section~3.2.2, Algorithm~A]{Knuth97}) in the Concorde TSP Solver with the Mersenne Twister \citep{Nishi04}, which is specifically designed for 64-bit machines and has a much longer period of $2^{19937}-1$ (compared to the default generator's period of $2^{55} - 1$). This modification was necessary, as our experiment to verify the curse of scale-freeness required running the RMS and ILS algorithms for an extremely long time.

We obtained the results presented in Figures~\ref{fig_RMS_1000cities} and \ref{fig_GSR_RMS_1000cities} by running $10^6$ trials for each of the six RMS algorithms (see Table~\ref{tb_RMS_algorithms}) on 100 random TSP instances (introduced in Section~\ref{subsec_instances}). For each set of $10^6$ trials, a unique random seed was assigned to the Mersenne Twister, and the Concorde TSP Solver used the corresponding sequence of random numbers.

We also obtained the results shown in Figures~\ref{RMS-ILS_RA-LKH}, \ref{RMS-ILS_NN-LKH}, and \ref{RMS-ILS_GR-LKH} by applying a similar procedure to each of the five TSPLIB instances listed in Table~\ref{TSPLIB_instances}.

\begin{enumerate} 
\item 100 runs of $10^6$ trials for each of the three RMS algorithms (RA + LK, NN + LK, and GR + LK); and
\item 100 runs of ``one initial start and $(10^6 - 1)$ random-walk kicks'' for each of the three ILS algorithms (RA + LK, NN + LK, and GR + LK). 
\end{enumerate}
The introduction of randomness into these RMS and ILS algorithms followed the same methodology used to generate Figures~\ref{fig_RMS_1000cities} and \ref{fig_GSR_RMS_1000cities}.

\if0
%%%%%%%%%%%%%%%%%%%%%%%%%% Acknowledgments %%%%%%%%%%%%%%%%%%%%%%%%%%

\section*{Acknowledgments}
The research of Hiroyuki Masuyama was supported in part by JSPS KAKENHI Grant Number JP21K11770.
\fi
%%%%%%%%%%%%%%%%%%%%%%%%%%%%%%%%%%%%%%%%%%%%%%%%%%%%%%%%%%%%%%%%%%%%%%
%%%							Springer
%%\bibliographystyle{sn-mathphys-ay} 
%%%%%%%%%%%%%%%%%%%%%%%%%%%%%%%%%%%%%%%%%%%%%%%%%%%%%%%%%%%%%%%%%%%%%%
%%%							Elsevier
%
\bibliographystyle{elsarticle-harv}
%%%%%%%%%%%%%%%%%%%%%%%%%%%%%%%%%%%%%%%%%%%%%%%%%%%%%%%%%%%%%%%%%%%%%%
%%%							INFORMS
%\bibliographystyle{informs2014}
%%%							Plain
%\bibliographystyle{plain} % plain, alpha, abbrv, unsrt
%%%%%%%%%%%%%%%%%%%%%%%%%%%%%%%%%%%%%%%%%%%%%%%%%%%%%%%%%%%%%%%%%%%%%%
\bibliography{OPT-EVT}

\end{document}